\RequirePackage{etoolbox}
\csdef{input@path}{{style/}{graphics/}}
\documentclass[aos,MSNbibl,seceqn,dvips]{arximspdf}
\usepackage{mathbh,dcolumn}
\usepackage{graphicx}

\doi{10.1214/14-AOS1291}
\volume{43}
\issue{2}
\pubyear{2015}
\firstpage{519}
\lastpage{545}
\docsubty{FLA}

\makeatletter
\newcolumntype{d}[1]{D{.}{.}{#1}}
\newcommand{\rrVert}{\Vert}
\newcommand{\llVert}{\Vert}
\def\cal{\mathcal}

\newtheorem{thmm}{Theorem}[section]
\newtheorem{prop}[thmm]{Proposition}

\newproclaim{remark}{Remark}[section]
\newproclaim{example}{Example}
\newcommand{\eqref}[1]{(\ref{#1})}
\newcommand{\bZ}{\bar{Z}}
\newcommand{\C}{{\cal C}}
\newcommand{\D}{{\cal D}}
\newcommand{\F}{{\cal F}}
\newcommand{\N}{{\cal N}}
\newcommand{\R}{{\cal R}}

\newcommand{\chga}{\check{\gamma}}

\newcommand{\al}{\alpha}
\newcommand{\be}{\beta}
\newcommand{\de}{\delta}
\newcommand{\ep}{\varepsilon}
\newcommand{\Ga}{\Gamma}
\newcommand{\la}{\lambda}
\newcommand{\ka}{\kappa}
\newcommand{\om}{\omega}
\newcommand{\Si}{\Sigma}
\newcommand{\te}{\theta}
\newcommand{\Th}{\Theta}

\newcommand{\hga}{\hat{\gamma}}
\newcommand{\hgan}{{\hat{\gamma}}_n}
\newcommand{\hsi}{\hat{\sigma}}
\newcommand{\hsin}{\hat{\sigma}_n}
\newcommand{\hthn}{{\hat{\theta}}_{n}}
\newcommand{\Var}{\operatorname{Var}}

\newcommand{\ci}{{\iota}}

\newcommand{\ind}{\mathbh{1}}

\newcommand{\nb}{\nabla}

\newcommand{\nti}{n\rightarrow\infty}
\newcommand{\raw}{\mathop\rightarrow}

\newcommand{\cid}{\stackrel{d}{\longrightarrow}}

\newcommand{\tG}{\tilde{G}}
\newcommand{\tI}{\tilde{I}}

\renewcommand{\rho}{\sigma}
\newcommand{\vr}{\varrho}
\def\be{\begin{eqnarray*}
} \def\ee{
\end{eqnarray*}}
\def\L{\mathcal{L}}
\def\a{\alpha}
\def\b{\beta}
\def\g{\gamma}
\def\D{\mathcal{D}}
\def\R{\mathcal{R}}
\def\x{{\mathbf{x}}}

\def\w{\bolds{\omega}}

\def\s{{\mathbf{s}}}
\def\h{{\mathbf{h}}}
\def\0{{\mathbf{0}}}
\def\h{{\mathbf{h}}}
\def\F{\mathcal{F}}
\def\cos{\operatorname{cos}}
\newcommand{\co}{c_*}
\newcommand{\px}{P_{\mathbf{X}}}

\def\bbn{\mathbb{N}}
\def\bbr{\mathbb{R}}
\def\bbz{\mathbb{Z}}

\def\bbi{\mathbb{I}}

\def\bbX{\mathbb{X}}

\makeatother

\begin{document}
\begin{frontmatter}

\title{A frequency domain empirical likelihood method for irregularly
spaced spatial data}
\runtitle{Spatial empirical likelihood}

\begin{aug}
\author[A]{\fnms{Soutir} \snm{Bandyopadhyay}\corref{}\thanksref{T1}\ead[label=e1]{sob210@lehigh.edu}},
\author[B]{\fnms{Soumendra N.} \snm{Lahiri}\thanksref{T2}\ead[label=e2]{snlahiri@ncsu.edu}}\\
\and
\author[C]{\fnms{Daniel J.} \snm{Nordman}\thanksref{T3}\ead[label=e3]{dnordman@iastate.edu}}
\thankstext{T1}{Supported in part by NSF Grant DMS-14-06622.}
\thankstext{T2}{Supported in part by NSF Grants DMS-10-07703 and DMS-13-10068.}
\thankstext{T3}{Supported in part by NSF Grant DMS-14-06747.}
\runauthor{S. Bandyopadhyay, S. N. Lahiri and D. J. Nordman}
\affiliation{Lehigh University, North Carolina State University and
Iowa~State~University}
\address[A]{S. Bandyopadhyay\\
Department of Mathematics\\
Lehigh University\\
Bethlehem, Pennsylvania 18015\\
USA\\
\printead{e1}}
\address[B]{S. N. Lahiri\\
Department of Statistics\\
North Carolina State University\\
Raleigh, North Carolina 27695-8203\\
USA\\
\printead{e2}}
\address[C]{D. J. Nordman\\
Department of Statistics\\
Iowa State University\\
Ames, Iowa 50011\\
USA\\
\printead{e3}}
\end{aug}

%
\received{\smonth{12} \syear{2013}}
%
\revised{\smonth{11} \syear{2014}}

%
\begin{abstract}
This paper develops empirical likelihood methodology for irregularly
spaced spatial data in the frequency domain. Unlike the frequency
domain empirical likelihood (FDEL) methodology for time series (on a
regular grid), the formulation of the spatial FDEL needs special care
due to lack of the usual orthogonality properties of the discrete
Fourier transform for irregularly spaced data and due to presence of
nontrivial bias in the periodogram under different spatial asymptotic
structures. A spatial FDEL is formulated in the paper taking into
account the effects of these factors. The main results of the paper
show that Wilks' phenomenon holds for a scaled version of the logarithm
of the proposed empirical likelihood ratio statistic in the sense that
it is asymptotically distribution-free and has a chi-squared limit. As
a result, the proposed spatial FDEL method can be used to build
nonparametric, asymptotically correct confidence regions and tests for
covariance parameters that are defined through spectral estimating
equations, for irregularly spaced spatial data. In comparison to the
more common studentization approach, a~major advantage of our method is
that it does not require explicit estimation of the standard error of
an estimator, which is itself a very difficult problem as the
asymptotic variances of many common estimators depend on intricate
interactions among several population quantities, including the
spectral density of the spatial process, the spatial sampling density
and the spatial asymptotic structure. Results from a numerical
study are also reported to illustrate the methodology and its finite sample
properties.
\end{abstract}

%

\end{frontmatter}

\section{Introduction}\label{sec1}
\label{intro}
In recent years, there has been a surge in research interest in the
analysis of spatial data using the frequency domain approach; see, for
example, Hall and Patil \cite{hallpatil1994}, Im, Stein and Zhu
\cite
{im2007}, Fuentes \cite{fuentes2006,fuentes2007}, Matsuda and Yajima
\cite{matsuda2009} and the
references therein.
An intent of frequency domain analysis is to allow for inference about
covariance structures through a data transformation and possibly
without a full spatial model, though this approach has complications.
In
contrast to the time series case where observations are usually taken
at regular points in time, the data sites are typically irregularly
spaced for random processes observed over space. The lack of a fixed
spacing and possible nonuniformity of the (irregularly spaced)
data-locations destroy the orthogonality properties of the sine- and
cosine-transforms of the data, making Fourier analysis in such
problems a challenging task. In a recent paper, Bandyopadhyay and
Lahiri \cite{bandy2010}
(hereafter referred to as [BL]) carried out a detailed investigation
of the properties of a suitably defined discrete Fourier transform
(DFT) of irregularly spaced spatial data, and provided a
characterization of the asymptotic independence property of
the spatial DFTs. In this paper, we utilize the insights and
findings of [BL] to formulate a frequency domain empirical
likelihood
(FDEL) for such spatial data. The FDEL method is shown to admit a
version of
the Wilks' theorem for test statistics about spatial covariance
parameters (e.g., having chi-square limits similarly to parametric
likelihood), without explicit assumptions on the data distribution or
the spatial sampling design.

To highlight potential advantages of the FDEL approach in this context,
suppose that
$\{Z(\mathbf{s})\dvtx\mathbf{s}\in\bbr^d\}$ ($d\in\bbn\equiv\{
1,2,\ldots\}
$) is a
zero mean
second-order stationary process that is observed at (irregularly spaced)
locations $\mathbf{s}_1,\ldots, \mathbf{s}_n$ in a domain $\D
_n\subset\bbr^d$.
Also, suppose that we are interested in fitting a parametric variogram
model $\{\check{\gamma}(\cdot;\te) \dvtx\te\in\Th\}$, $\Th\in
\bbr^p$
($p\in\bbn$) using
the least squares approach (cf. Cressie \cite{cressie1993}). A spatial
domain approach
is based on estimating the parameter $\te$ using
\[
\tilde{\te}_n = \operatorname{argmin} \Biggl\{ \sum
_{i=1}^m \bigl(2\tilde{\gamma}_n(
\mathbf{h}_i) - \check{\gamma}(\mathbf{h}_i;\te)
\bigr)^2 \dvtx\te\in\Th \Biggr\},
\]
where $\mathbf{h}_1,\ldots, \mathbf{h}_m$ are some user specified
lags and where
$2 \tilde{\gamma}_n(\mathbf{h}_i)$ is a nonparametric estimator of the
variogram of the $Z(\cdot)$-process at lag $\mathbf{h}_i$. Since the
data locations
$\mathbf{s}_1,\ldots, \mathbf{s}_n$ are irregularly spaced, a
nonparametric estimator
$2\tilde{\gamma}_n(\cdot)$ of the variogram typically requires
\textit{smoothing} which results in
a slow rate of convergence, particularly in dimensions $d\geq2$.
Further, the asymptotic variance of $\tilde{\te}_n$ in such situations
involves the spectral density
of $Z(\cdot)$-process
and the spatial sampling density of the data-locations $\mathbf
{s}_1,\ldots,
\mathbf{s}_n$
(cf. Lahiri and Mukherjee \cite{lahiri2004}) which must be estimated
from the data
to carry out inference on $\te$ using the asymptotic distribution.
In contrast, the FDEL approach completely bypasses the need to estimate
$2\check{\gamma}(\cdot)$ directly and it also carries out an automatic
adjustment
for the complicated asymptotic variance term
in its inner mechanics, producing a distribution-free limit law that
can be readily used for constructing valid tests and confidence regions
for $\te$. See Example~\ref{ex3} in Section~\ref{fdel}
for more details of the FDEL construction in this case and Section~\ref{sec6.2}
for a data example demonstrating the advantages of the proposed
spatial FDEL method over the traditional spatial domain
approach. In general, the proposed FDEL method provides
a nonparametric ``likelihood''-based inference method
for covariance parameters of a spatial process observed at
irregularly spaced spatial data-locations without requiring
specification of a parametric joint data model.

Originally proposed by Owen \cite{owen1988,owen1990} for independent
observations, empirical likelihood (EL) allows for nonparametric
likelihood-based inference in a broad range of applications (Owen
\cite
{owen2010}),
such as construction of confidence regions for parameters
that may be calibrated through the asymptotic chi-squared distribution
of the log-likelihood ratio. This is commonly referred to as the
\textit{Wilks' phenomenon}, in analogy to the asymptotic distributional
properties of likelihood ratio tests in traditional
parametric problems (Wilks \cite{wilks1938}). In particular, EL does
not require any direct estimation of variance or skewness (Hall and
La Scala \cite{hall1990}). However, a~difficulty with extending EL
methods to dependent data is then to ensure that ``correct'' variance
estimation occurs automatically within the mechanics of EL under
dependence. For (regularly spaced) time series data,
this is often accomplished by using a blockwise
empirical likelihood (BEL) method (cf. Kitamura~\cite{kitamura1997}),
which was further extended to the case of spatial data observed on a
regular grid by Nordman \cite{nordman2008a,nordman2008b} and Nordman
and Caragea \cite{nordman2008c}.

Monti \cite{monti1997} and Nordman and Lahiri \cite{nordman2006}
proposed periodogram-based EL methods for time series data. Their works
show that, in view of the asymptotic independence of the DFTs,
an analog of the EL formulation for independent data satisfies
Wilks' phenomenon in the frequency domain. As a result, the vexing
issue of block length choice can be completely avoided by working
with the DFTs of (regularly spaced) time series data. In this paper,
we extend the frequency domain approach to irregularly spaced
spatial data. Such an extension presents a number of unique
challenges that are inherently associated with the spatial
framework.
First, the irregular spacings of the data locations make the
usefulness of the DFT itself questionable, as the basic
orthogonality property of the sine- and cosine-transforms of
gridded data at Fourier frequencies [i.e., at frequencies
$\om_j=2\pi j/n$ for $j=0,1,\ldots,(n-1)$ for
a time series sample of size $n$] no longer holds
(cf.~[BL]).
Second,
unlike the \textit{compact} frequency domain $[0,2\pi]$ for regular
time series, in the case of irregularly spaced spatial processes
sampled in $d$-dimensional Euclidean space $\bbr^d$, one must deal
with the \textit{unbounded} frequency domain $\bbr^d$.
Third,
as noted in Matsuda and Yajima \cite{matsuda2009} and [BL], the
periodogram of
irregularly spaced spatial data can be severely biased
(for the spectral density)
and must be pre-processed.
Finally,
in contrast to the unidirectional flow of time that drives the
asymptotics in
the time series case, for irregularly spaced spatial data on an increasing
domain, more than one possible asymptotic structure can arise
depending on
the relative growth rates of the volume of the sampling region and
the sample size (cf. Cressie \cite{cressie1993}, Hall and Patil
\cite
{hallpatil1994}, Lahiri \cite{lahiri2003a}).
A desirable property of any FDEL method for irregularly spaced
spatial data would be to guarantee Wilks' phenomenon for the spatial
FDEL ratio statistic with minimal or no explicit adjustments for
the different asymptotic regimes. This would ensure a sort of robustness
property for the spatial FDEL and would allow the user to use the method
in practice without
having to explicitly tune it for the effects of different spatial
asymptotic structures, which is often not very obvious for a given
data set at hand (cf. Zhang and Zimmerman \cite{zhang2005}).

To motivate the construction of our spatial FDEL (hereafter SFDEL),
first we briefly review some relevant results (cf. Section~\ref{prelim})
that provide crucial insights into the properties of the DFT and
periodogram of irregularly spaced spatial data under different
spatial asymptotic structures.
Our main result is the asymptotic chi-squared distribution of the
SFDEL ratio statistic under fairly general regularity conditions on the
underlying spatial process. However, it turns out that the spatial
asymptotic structure has a nontrivial and nonstandard effect on the
limit law. When the spatial sample size $n$ grows at a rate comparable
to the volume of the sampling region, we shall call this the \textit
{pure increasing domain} or PID asymptotic structure, while a faster
growth rate of $n$ (due to infilling) will be called the \textit{mixed
increasing domain} or MID asymptotic structure (see Section~\ref
{prelim} for more details). To describe the peculiarity of the limit
behavior of the SFDEL, let $\R_n(\te_0)$ denote the SFDEL ratio
statistic for a covariance parameter of interest $\te\in\bbr^p$ under
$H_0\dvtx\te=\te_0$ based on a sample of size $n$.
The main results of the paper show that under some regularity conditions,
%
%
\begin{equation}
- 2\log\R_n(\te_0) \stackrel{d} {\rightarrow}
\chi^2_p \label{MID-f}
\end{equation}
under MID with a sufficiently fast rate of infilling. In contrast,
under PID
and under MID with a relatively slow rate of infilling,
one gets
%
%
\begin{equation}
- 2\log\R_n(\te_0) \raw^d 2
\chi^2_p. \label{PID}
\end{equation}
Thus, the limit distribution of $- 2\log\R_n(\te_0)$ here
changes from the more familiar $\chi^2_p$ to
a nonstandard $2\chi^2_p$ distribution, which
points to the intricacies associated with spatial asymptotics.
The main reason behind this strange behavior of the SFDEL ratio
statistic is the \textit{differential growth rates of two components}
in the variance term of the DFTs of irregularly spaced spatial data,
which alternate in their roles as the dominating term
depending on the strength of the infill component.

To overcome the dichotomous limit behavior of $- 2\log\R_n(\te_0)$
in \eqref{MID-f} and~\eqref{PID},
we construct a data based scaling $a_n=a_n(\te_0)$ (say) and show
that the rescaled version,
$-a_n 2\log\R_n(\te_0)$ attains the \textit{same} $\chi^2_p$ limit,
irrespective of the underlying spatial asymptotic structure.
This provides a unified method for EL based inference on covariance
parameters for irregularly spaced spatial data.
In addition, the proposed SFDEL method accomplishes
two major goals of the EL method of Owen \cite{owen1988,owen1990}
for independent data:
\begin{longlist}[(ii)]
\item[(i)] it shares the strength of EL methods to incorporate
\textit{automatic variance estimation} for spectral parameter
inference in its mechanics under different spatial asymptotic structures
and, at the same time,
\item[(ii)] it \textit{avoids} the difficult issue of block length selection.
\end{longlist}
A direct solution to either of these problems (i.e., explicit
variance estimation and optimal block length selection)
in the spatial domain is utterly difficult due to
highly complex effects of the
irregular spacings of the data sites and the spatial asymptotic
structures (cf. Lahiri \cite{lahiri2003a}, Lahiri and Mukherjee
\cite
{lahiri2004}) and
due to potentially nonstandard shapes of the sampling regions
(Nordman and Lahiri \cite{nordman2004}, Nordman, Lahiri and
Fridley \cite{nordman2007}).
Results from a simulation study in Section~\ref{sim}
show that accuracy of
the SFDEL method with the data-based rescaling is
very good even in moderate samples.

The rest of the paper is organized as follows. In Section~\ref{prelim},
we describe the theoretical framework and some preliminary
results on the properties of the DFT for
irregularly spaced spatial data from [BL] that play a
crucial role in the formulation of the SFDEL method. We
describe the SFDEL method in Section~\ref{fdel} and
give some examples of useful spectral estimating equations.
We state the regularity conditions and the main results
of the paper in Sections \ref{reg-cond} and Section~\ref{asy},
respectively.
Results from a simulation study and an illustrative data example
are given in Section~\ref{sim}. Proofs of the main results are
presented in Section~\ref{proof_fdel}.
Further details of the proofs and some additional simulation results are
given in the supplementary material~\cite{supp}.

%
\section{Preliminaries}
\label{prelim}
\subsection{Spatial sampling design}
\label{design}
Suppose that for each $n\geq1$ (where $n$ denotes the sample size),
the spatial process $Z(\cdot)$ is observed at data
locations $\mathbf{s}_1,\ldots,\mathbf{s}_n$ over a sampling region
$\D
_{n}\subset
\bbr^d$. We shall suppose that $\D_n$ is obtained by inflating a
prototype set $\D_{0}$ by a scaling factor $\lambda_n \in[1,\infty
)$ as
%
%
\begin{equation}
\D_{n}=\lambda_{n}\D_{0}, \qquad n\geq1,
\label{samp-r}
\end{equation}
where (as the most relevant prototypical case) $\D_{0}$ is an open
connected subset of $(-1/2,1/2]^d$ containing the origin and where
$\lambda
_{n}\uparrow\infty$ as $n\rightarrow\infty$ with $\lambda_n \gg
n^\varepsilon$ for some $\varepsilon>0$. Note that this is a common
formulation, allowing the sampling region $\D_n$ to have a variety of
shapes, such as polygonal, ellipsoidal and star-shaped regions that
can be nonconvex. In practice, $\lambda_n$ can be determined by the
diameter of a sampling region for use here (cf. Garc\'{i}a-Soid{\'a}n
\cite{gsoidan2007}, Hall and Patil \cite{hallpatil1994},
Maity and Sherman \cite{maity2012},
Matsuda and Yajima \cite{matsuda2009}). Let $\bbz=\{0,\pm1,\pm
2,\ldots\}$. To avoid pathological cases, we require that for any
sequence of real numbers $\{a_{n}\}_{n\geq1}$ such that
$a_{n}\rightarrow0+$ as $n\rightarrow\infty$, the number of\vspace*{1pt} cubes of
the form $a_{n}(\mathbf{j}+[0,1)^{d}), \mathbf{j}\in\mathbb{Z}^{d}$
that intersect
both $\D_{0}$ and $\D_{0}^{c}$ is of the order $O([a_{n}]^{-(d-1)})$ as
$n\rightarrow\infty$. This boundary condition holds for most regions
of practical interest.
We also suppose that
the irregularly spaced
data locations $\mathbf{s}_1,\ldots, \mathbf{s}_n\in\D_n$
are generated by a stochastic sampling design,
as
\[
\s_{i}\equiv\s_{in} = \lambda_{n}
\mathbf{X}_{i},\qquad1\leq i\leq n,
\]
where
$ \{{\mathbf{X}}_{k} \}_{k\geq1}$ is a sequence of
independent and identically distributed (i.i.d.) random vectors
with probability density $f(\x)$ with support $\mbox{cl.}(\D_0)$,
the {closure} of~$\D_{0}$.
Note that this formulation allows the number of sampling sites
to grow at a \textit{different rate} than the volume of the
sampling region, leading to different
asymptotic structures (cf. Cressie \cite{cressie1993}, Lahiri \cite
{lahiri2003a}).
When $n/\lambda_{n}^{d}\rightarrow\co\in(0,\infty)$, one gets
the PID asymptotic structure while for
$n/\lambda_{n}^{d}\rightarrow\infty$ as $n\rightarrow\infty$,
one gets the MID asymptotic structure.
Limit laws of common estimators are known to
depend on the spatial asymptotic structure; see Cressie \cite{cressie1993},
Du, Zhang and Mandrekar~\cite{du2009}, Lahiri and Mukherjee~\cite{lahiri2004}, Loh \cite{loh2005}, Stein \cite{stein1989b}
and the references therein.
%

\subsection{Spatial periodogram and its properties}
\label{prop-dft}
Define the DFT $d_n(\w)$ and the periodogram $I_n(\w)$ of $\{
Z(\mathbf{s}
_1),\ldots, Z(\mathbf{s}_n)\}$ at $\w\in\mathbb{R}^{d}$ as
%
%
\begin{equation}
{d}_{n}(\w)= \lambda_{n}^{d/2}n^{-1}
\sum_{j=1}^{n}Z(\s_{j}) \exp
\bigl(\ci\w{'}\s_{j} \bigr) \quad\mbox{and}\quad
{I}_{n}(\w)=\bigl |{d}_{n}(\w) \bigr|^{2}, \label{dft-p}
\end{equation}
where $\ci=\sqrt{-1}$.
In an equi-spaced time series, formulation and properties of
FDEL critically depend on the asymptotic independence
of the DFTs (cf. Brockwell and Davis \cite{brockwell1991}, Lahiri
\cite{lahiri2003b})
at the Fourier frequencies:
$w_j = 2\pi j/n, j=1,\ldots,n$ where $n$ is the sample size.
In a recent paper,
[BL] showed that the spatial DFTs [in \eqref{dft-p}] at two
sequences of frequencies $\{\bolds{\omega}_{1n}\}_{n\geq1}$,
$\{\bolds{\omega}_{2n}\}_{n\geq1} \subset{\mathbb{R}^{d}}$ are
\textit
{asymptotically
independent} (i.e., the joint limit law is a product of marginal
limits) \textit{if and only if the frequency sequences are asymptotically
distant}:
%
%
\begin{equation}
\bigl\|\lambda_n(\bolds{\omega}_{1n} - \bolds{
\omega}_{2n}) \bigr\| \raw \infty\qquad\mbox{as }\nti. \label{a-dist}
\end{equation}

This suggests that in analogy to the time series FDEL (i.e., using that
DFTs are approximately independent so that the independent data version
of EL may be applied to resulting periodogram values),
the formulation of spatial FDEL should preferably be
based on DFTs at a collection of frequencies that are well-separated.
A second important finding in [BL] is that unlike the case of
the equi-spaced time series data, the spatial periodogram $I_n(\cdot)$
has a nontrivial bias, depending on the spatial asymptotic structure.
In particular,
[BL] shows that
\[
EI_n(\bolds{\omega}) = \bigl[n^{-1}\lambda_n^d
{\sigma}(\0) + K \phi(\w) \bigr] \bigl(1+o(1) \bigr)
\]
for all $\bolds{\omega}\in{\mathbb{R}^{d}}$, where ${\sigma}(\cdot
)$ and
$\phi(\cdot)$
are respectively the autocovariance and the spectral density
functions of the $Z(\cdot)$-process
and where $K= (2\pi)^d \times\break \int_{{\mathbb{R}^{d}}} f^2 (\bolds{\omega
}) \,d\bolds{\omega}
$. As a result,
the spatial periodogram $I_n(\cdot)$ has a nontrivial bias [for
estimating $K\phi(\cdot)$]
at \textit{all} frequencies under PID, while the bias vanishes \textit
{asymptotically}
under MID. However, the quality of estimation of
the spectral density (up to the scaling by $K$)
improves under both PID and MID through an
explicit bias correction. Accordingly, we define the bias
corrected periodogram
%
%
\begin{equation}
\tilde{I}_{n}(\w)={I}_{n}(\w) - n^{-1}
\lambda_n^d \hsin(\0), \qquad\w\in\bbr^d,
\end{equation}
where $\hsin(\0) = n^{-1}\sum_{i=1}^n (Z(\mathbf{s}_i) - \bar
{Z}_n)^2$ is the
sample variance, with $\bar{Z}_n = n^{-1}\sum_{i=1}^n Z(\mathbf{s}_i)$
denoting the sample mean. We shall use $\tilde{I}_n(\cdot)$ in our
formulation of the SFDEL in the next section.

\section{The SFDEL method}
\label{fdel}
\subsection{Description of the method}\label{fdel.1}
For i.i.d. random variables, Qin and Lawless~\cite{qin1994} extended
the scope
of Owen's \cite{owen1988} original formulation, linking estimating
equations and EL, and developed EL methodology for
such parameters. In a recent work, Nordman and Lahiri \cite
{nordman2006} (hereafter
referred to as [NL]) formulated a FDEL for inference on parameters
of an equi-spaced time series defined through
spectral estimating equations (i.e., estimating equations in
the frequency domain $[-\pi,\pi]$). In a similar spirit,
we now define the SFDEL for parameters $\theta\in\Th\subset\bbr^p$,
defined through spectral estimating equations (but now defined
over all of $\bbr^d$).
Specifically, let $G\dvtx\bbr^d\times\Th\rightarrow\mathbb{R}^{p}$
denote a vector of bounded estimating functions such that
$G_{\theta}(\cdot) \equiv G(\cdot;\theta)$ satisfies the spectral
moment condition
%
%
\begin{equation}
\int_{{\mathbb{R}^{d}}} G_{\theta}(\bolds{\omega}) \phi(\bolds {
\omega}) \,d\bolds{\omega}=0, \label{spec-mc}
\end{equation}
where recall that $\phi(\cdot)$ denotes the spectral density of the
process $Z(\cdot)$.
Because of their use in the SFDEL method to follow [cf. (\ref{Fdel})],
we refer to the functions
$G_{\theta}(\bolds{\omega})$ as estimating functions, though these
are not
functions of data directly but
rather of parameters $\theta\in\Th$ and frequencies $\bolds{\omega
}\in
\mathbb{R}^d$.
In view of the symmetry of the
spectral density $\phi(\cdot)$, without loss of generality (w.l.g.),
we shall assume that
$G_{\theta}(\cdot)$ is symmetric about zero, that is,
$G_{\theta}(\w)= G_{\theta}(-\w)$ for all $\w\in\bbr^d$.
An asymmetric $G_{\theta}(\cdot)$ can always be symmetrized, as in
Example~\ref{ex2}
of Section~\ref{examples} below where we give examples of $G_{\theta
}(\cdot)$
in some important inference problems. 

The SFDEL defines a nonparametric likelihood for the
parameter $\te$ using a \textit{discretized} sample version of
the above spectral moment condition.
Accordingly, for $\kappa\in(0,1)$, $\eta\in[\ka,\infty)$ and
$C^*\in(0,\infty)$, let
%
%
\begin{equation}
\N=\N_{n}= \bigl\{\mathbf{j}\lambda_{n}^{-\kappa}
\dvtx\mathbf{j} \in\mathbb{Z}^{d}, \mathbf{j}\in \bigl[-C^*
\lambda_{n}^{\eta}, C^*\lambda_{n}^{\eta}
\bigr]^{d} \bigr\} \label{grid}
\end{equation}
be the set of discrete frequencies, where $\lambda_n$ is as in
\eqref{samp-r}. Let $N=|\N|$ be the size of $\N$. For notational
convenience, also denote the elements of $\N$ by $\w_{kn},\
k=1,\ldots,N$ (with an arbitrary ordering of
the $N$ elements of $\N$). The frequency grid has two important
qualities. First, since $\ka<1$, for any $j\neq k$,
the sequences $\{\bolds{\omega}_{jn}\}$ and $\{\bolds{\omega}_{kn}\}
$ are
asymptotically distant [cf. \eqref{a-dist}], guaranteeing their
associated periodogram values are approximately independent.
Further, $\{\w_{1n},\ldots,\w_{Nn}\}$ forms a \textit{regular}
lattice over the hyper-cube $ [-C^*\lambda_{n}^{\eta-\kappa},
C^*\lambda_{n}^{\eta-\kappa} ]^{d}$, with spacings of
length $\lambda^{-\kappa}_{n}$ in each direction, and
$[-C^*\lambda_n^{\eta-\ka}, -C^*\lambda_n^{\eta-\ka}]^d
\uparrow{\mathbb{R}^{d}}$ as $\lambda
_n\uparrow\infty$ when
$\nti$ for $\eta>\ka$, covering the entire range of
the integral in (\ref{spec-mc}) in the limit. That is, the frequency
grid expands to necessarily cover the entire frequency domain ${\mathbb
{R}^{d}}$ of
interest.
The exact conditions on $\ka$ and $\eta$ are specified
in Section~\ref{reg-cond} below.

Now using the frequencies $\{\w_{kn},\ k=1,\ldots,N\}$,
we define the SFDEL function for $\theta$ by
%
%
\begin{eqnarray}\label{Fdel}
&&\L_{n}(\theta)=\sup \Biggl\{\prod_{k=1}^{N}p_{k}
\dvtx\sum_{k=1}^{N}p_{k}=1,
p_{k}\geq0 \mbox{ and}
\nonumber
\\[-8pt]
\\[-8pt]
\nonumber
&& \hspace*{76pt}\sum_{k=1}^{N}p_{k}G_{\theta}(
\w_{kn}) \tilde{I}_{n}(\w_{kn})=0 \Biggr\},
\end{eqnarray}
provided that the set of $p_k$ satisfying the conditions on the right-hand
side is nonempty. When no such $\{p_k\}$ exists,
$\L_n (\theta)$ is defined to be $0$. We note that the computation
of (\ref
{Fdel}) is the same as in EL formulations for independent data; see
Owen \cite{owen1990,owen2010} and Qin and Lawless \cite{qin1994} for
these details.

Next, note that without the spectral moment constraint,
$\prod_{k=1}^{N}p_{k}$ attains its maximum when each $p_{k}=1/N$.
Hence, we define the SFDEL ratio
statistic for testing the hypothesis $H_0\dvtx\theta= \theta_0$ as
\[
\R_{n}(\theta_0)=\mathcal{L}_{n}(
\theta_0)/\bigl(N^{-N}\bigr).
\]
The SFDEL test rejects $H_0$ for small values of $\R_{n}(\theta_0)$.
Similarly, one can use the SFDEL method to construct confidence regions for
$\te$ using the large sample distribution of the SFDEL ratio
statistic. In Section~\ref{reg-cond}, we state a set of regularity
conditions that will
be used for deriving the limit distribution of $-2\log\R_n(\theta_0)$.
This, in particular, would allow one to calibrate the SFDEL tests
and confidence regions in large samples.

%
\subsection{Examples of estimating equations}
\label{examples}
We now give some examples of spectral estimating equations for
parameters of interest in frequency domain analysis (cf. Brockwell and
Davis \cite{brockwell1991},
Cressie \cite{cressie1993},
Journel and Huijbregts \cite{journel1978},
Lahiri, Lee and Cressie \cite{lahiri2002}).

\begin{example}[(Autocorrelation)]\label{ex1}
Suppose that
we are interested in nonparametric estimation of
the autocorrelation of the $Z(\cdot)$-process at lags
$\h_{1},\ldots,\h_{p}$ for some $p\geq1$.
Then $\te= (\vr(\mathbf{h}_1),\ldots,\vr(\mathbf{h}_p))'$
with
$
\vr(\mathbf{h})=\break \operatorname{corr}(Z(\h),Z(\0)) = \int\cos{(\h
^{\prime}\w)}
\phi(\w)\,d\w/\int\phi(\w)\,d\w
$
where $A'$ denotes the\break transpose of a matrix $A$. Thus, in this case,
%
%
\begin{equation}
G_{\theta}(\w)= \bigl(\cos{ \bigl(\h_{1}^{\prime}\w
\bigr)},\ldots,\cos{ \bigl(\h_{p}^{\prime}\w \bigr)}
\bigr)^{\prime}-\theta. \label{efn-cr}
\end{equation}
Estimating functions can also be formulated with hypothesized
autocorrelations (e.g., white noise)
to set-up goodness-of-fit tests in the SFDEL approach, in the spirit of
Portmanteau tests Brockwell and Davis \cite{brockwell1991}.
\end{example}

\begin{example}[(Spectral distribution function)]\label{ex2}
For $\mathbf{t}=(t_1,\ldots,t_d)'\in\mathbb{R}^{d}$, let
\[
\Phi^0(\mathbf{t})=\int\ind_{(-\infty,\mathbf{t}]}(\w)\phi(\w )\,d\w\Big/\int
\phi(\w)\,d \w
\]
denote the normalized spectral distribution function, where
$\ind(\cdot)$ denotes the indicator function and
$(-\infty,\mathbf{t}]=(-\infty,t_{1}]\times\cdots\times(-\infty,t_{d}]$.
The function $\Phi^0(\cdot)$ plays an important role in
determining the smoothness of the sample paths of the
random field $Z(\cdot)$ (cf. Stein \cite{stein1989b}).
Suppose that the parameter of interest is now given by
$\theta=(\Phi^0(\mathbf{t}_1),\ldots,\Phi^0(\mathbf{t}_p))'$ for
some given
set of vectors $\mathbf{t}_1,\ldots, \mathbf{t}_p \in{\mathbb
{R}^{d}}$. In
this case,
the relevant estimating function is $ G_{\theta}(\w) =
[\tG_{\theta}(\w)+\tG_{\theta}(-\w)]/2$, $\w\in\bbr^d$, where
%
%
\begin{equation}
\tG_{\theta}(\w)= \bigl(\ind_{(-\infty,\mathbf{t}_1]}(\w ),\ldots,
\ind_{(-\infty,\mathbf{t}_p]}(\w) \bigr)' - \theta. \label{efn-df}
\end{equation}
\end{example}

\begin{example}[(Variogram model fitting)]\label{ex3}
A~popular approach
to fitting a parametric variogram model to spatial data is through
the method of least squares (cf. Cressie \cite{cressie1993}). Let
$\{2\chga(\cdot;\te) \dvtx\te\in\Th\}$, $\Th\subset\bbr^p$ be
a class of valid variogram models for the true variogram
$ 2\chga(\mathbf{h}) \equiv\Var(Z(\mathbf{h}) - Z(\0))$, $\mathbf
{h}\in{\mathbb{R}^{d}}$
of the spatial process. Let $2\gamma(\cdot;\te)
\equiv2\chga(\cdot;\te)/{\sigma}(\0)$ and
$2\gamma(\cdot)\equiv2\chga(\cdot)/{\sigma}(\0)$ denote
their scale-invariant versions, where
${\sigma}(\0)=\operatorname{Var}(Z(\0))$. Also, let $2\hga_n(\mathbf
{h})$ denote
the sample variogram at lag $\mathbf{h}$ based on
$Z(\mathbf{s}_1),\ldots,Z(\mathbf{s}_n)$ (cf. Chapter~2,
Cressie \cite{cressie1993}), scaled by $\hsi_n(\0) =
n^{-1}\sum_{i=1}^n (Z(\mathbf{s}_i) -\bZ_n)^2$ where
$\bZ_n=n^{-1}\sum_{i=1}^n Z(\mathbf{s}_i)$. Then one can
fit the variogram model by estimating the
parameter $\te$ by
\[
\hthn= \operatorname{argmin} \Biggl\{ \sum_{i=1}^m
\bigl( 2\hgan(\mathbf{h}_i) - 2 \gamma(\mathbf{h}_i;
\te) \bigr)^2 \dvtx\te\in\Th \Biggr\}
\]
for a given set of lags $\mathbf{h}_1,\ldots, \mathbf{h}_m$. This
corresponds to minimizing the population
criterion $\sum_{i=1}^m ( 2\g(\mathbf{h}_i) - 2 \gamma(\mathbf
{h}_i;\te
) )^2$
which, under some
mild conditions, determines the true parameter $\te_0$ uniquely
(cf. Lahiri, Lee and Cressie \cite{lahiri2002}).
Under these conditions, $\te=\te_0$ is the unique solution to the equation
\[
\sum_{i=1}^m \bigl(2\gamma(
\mathbf{h}_i) - 2\gamma(\mathbf{h}_i;\te) \bigr) \nb
\bigl[2\gamma( \mathbf{h}_i;\te) \bigr] =0,
\]
where $\nb[2\gamma(\mathbf{h};\te)]$ denotes the $p\times1$ vector
of first-order partial derivatives of $2\gamma(\mathbf{h};\te)$
with respect to $\te$. Hence, expressing the variogram
in terms of the spectral density function, we get the
following equivalent spectral estimating equation:
%
%
\begin{equation}
\int \Biggl[\sum_{i=1}^m \bigl\{1 - \cos
\bigl(\mathbf{h}_i'\bolds{\omega} \bigr)-\gamma(
\mathbf{h};\te) \bigr\} \nb \bigl[2\gamma(\mathbf{h}_i;\te) \bigr]
\Biggr]\phi(\bolds{\omega }) \,d\bolds{\omega}=0, \label{vrg-eeq}
\end{equation}
which can be used for defining the SFDEL for $\te$.
As pointed out in Section~\ref{sec1}, the spatial domain approach yields
asymptotically correct confidence regions for $\te$
through asymptotic normal distribution of $\hthn$,
but it necessarily requires one to estimate
the limiting asymptotic variance and
is subject to the curse of dimensionality, resulting from
nonparametric smoothing in $d$-dimensions. In comparison,
the SFDEL can be applied with the spectral estimating
equation (\ref{vrg-eeq}) to produce asymptotically
correct confidence region for $\te$, without
explicit estimation of the standard error.
\end{example}

Note that the spectral estimating equation approach
can also be extended to estimation of $\te$ based on the
weighted- and the generalized-least squares
criteria (cf. Cressie \cite{cressie1993},
Lahiri, Lee and Cressie \cite{lahiri2002}), where in addition to the partial
derivatives, suitable weight matrices enter into the
corresponding versions of \eqref{vrg-eeq}. A similar advantage
of the spatial FDEL method continues to hold in these
cases.

In the next section, we introduce some notation and the regularity conditions
to be used in the rest of the paper.

\section{Regularity conditions}
\label{reg-cond}
\subsection{Notation and lemmas}\label{sec4.1}
First, we introduce some notation. For $x, y\in\bbr$,
let $x_+=\max\{x,0\}$, ${\lfloor{x}\rfloor}=$ the floor function of $x$,
$x\wedge y = \min\{x,y\}$ and $x\vee y = \max\{x,y\}$.
Let $\bbi_k$ denote the identity matrix of order $k$ ($k\geq1$). For
two sequences $\{s_{n}\}$ and $\{t_{n}\}$ in $(0,\infty)$, we write
$s_{n}\sim t_{n}$ if $\lim_{n\rightarrow\infty}s_{n}/t_{n}=1$. For
${\mathbf x}=(x_1,\ldots,x_k)'\in\mathbb{R}^{k}$, let $\llVert
{\mathbf
x}\rrVert
_{1}=|x_1|+\cdots+|x_k|$ and $\|{\mathbf x}\|=(|x_1|^2+\cdots+|x_k|^2)^{1/2}$
respectively denote the $\ell^{1}$- and $\ell^{2}$-norms of ${\mathbf x}$.
Also, let $d_1(E_{1},E_{2})=\inf\{\llVert\x-\s\rrVert_{1}\dvtx\x
\in
E_{1}, \s\in E_{2}\}$, $E_1,E_2\subset\mathbb{R}^{k}$.
For $a,b\in(0,\infty)$, define the strong mixing coefficient of
$Z(\cdot)$ as $\a(a;b)=\sup\{|P(A_{1}\cap A_{2})-P(A_{1})P(A_{2})|
\dvtx A_{i}\in\F_{Z}(E_{i}), E_{i}\in\mathbb{C}_{b}, i=1,2,
d_1(E_{1},E_{2}) \geq a\}$ where $\F_{Z}(E)={\sigma}\langle Z(\s
)\dvtx\s\in
E\rangle$ and $\mathbb{C}_b$ is the collection of $d$-dimensional
rectangles with volume $b$ or less.

As indicated earlier, we suppose that the random
field $\{Z(\s)\dvtx\s\in\mathbb{R}^{d}\}$ is \textit{second-order
stationary}
(but not necessarily strictly stationary) with zero mean and
autocovariance function
$
\rho(\cdot)
$
and
spectral density function $\phi(\cdot)$. Also, recall that the
scaling sequence $\la_n$ is as in \eqref{samp-r} and that
$\kappa$, $\eta$ and $N$ are as in Section~\ref{fdel.1}, specifying
the SFDEL grid in the frequency domain. Further,
the constant $c_*\equiv\lim_{\nti} n/\la_n^d$
determines the spatial asymptotic structure
where $c_*\in(0,\infty)$ for PID and $c_*=\infty$ for MID.
Write $c_n = n/\la_n^d$, $I_n^*(\w) = I_n(\w) - c_n^{-1}{\sigma
}(\mathbf{0})$
and $A_n(\w) = c_n^{-1}{\sigma}(\0)+ K \phi(\w)$, $\w\in\bbr^d$,
where $K= (2\pi)^d \int f^2$.
Let $\Si_n=2 \sum_{k=1}^N G_{\theta_0}(\w_{kn})
G_{\theta_0}(\w_{kn})' A_n^2 (\w_{kn})$.
Also, let\vspace*{1pt} $G_{j,\te_0}$ denote the $j$th component of
$G_{\te_0}$. Write $b_n^2 = Nc_n^{-2}+ \la_n^{\ka d}$.
From Section~\ref{proof_fdel}, it follows that $b_n^2$
gives a unified representation for growth rate
of the self-normalizing factor in the SFDEL
ratio statistic under different asymptotic structures
considered in the paper.

\subsection{Conditions}\label{sec4.2}
We are now ready to state the regularity conditions.

%
\begin{enumerate}[(C.3)]
\item[(C.0)] The strong mixing coefficient satisfies
$\a(a,b)\leq\g_{1}(a)\g_{2}(b)$, for any $a,b\in(0,\infty)$, with
respect to
some left continuous nonincreasing
function $\g_1\dvtx(0,\infty) \raw[0,\infty)$ and some
right continuous nondecreasing function $\g_2\dvtx (0,\infty) \raw
(0,\infty)$.

\item[(C.1)] There exist $\de\in(0,1]$ such that
$\zeta_{4+\de}\equiv\sup\{(E|Z(\s)|^{4+\de})^{{1}/{(4+\de)}}
\dvtx\mathbf{s}\in\bbr^d\} <\infty$ and
$\sum_{k=1}^\infty k^{3d}[\gamma_1(k)]^{{\de}/{(4+\de)}}<\infty$.
\item[(C.2)] (i) The spatial sampling density $f(\cdot)$ is
everywhere positive on $\D_0$ and satisfies a Lipschitz condition:
There exists a $C_0\in(0,\infty)$ such that
\[
\bigl|f({\mathbf x}) - f(\mathbf{y}) \bigr| \leq C_0 \|{\mathbf x}-\mathbf{y}
\| \qquad\mbox{for all } {\mathbf x},\mathbf{y}\in\D_0.\vspace*{-6pt}
\]
\begin{enumerate}[(ii)]
\item[(ii)] There exist $C_1\in(0,\infty)$ and $a_0\in(d/2,d]$ such
that 
\begin{eqnarray}
\biggl|\int e^{\ci\bolds{\omega}^{\prime}{\mathbf x}}f({\mathbf x})\,d{\mathbf x}
 \biggr| + \biggl|\int
e^{\ci\bolds{\omega}^{\prime}{\mathbf x}}f^2({ \mathbf x})\,d{\mathbf x} \biggr| \leq
C_1\|\bolds{\omega}\|^{-a_0}\nonumber \\
\eqntext{\mbox{for all } \|\bolds{
\omega}\|>C_1.}
\end{eqnarray}
\end{enumerate}
\item[(C.3)] (i) For each $j=1,\ldots,p$, $G_{j,\te_0}(\cdot)$ is
bounded, symmetric, and almost everywhere continuous on $\bbr^d$ (with
respect to the Lebesgue measure), and $\int G_{\te_0}(\w)\phi(\w
)\,d\w
=0$;
\begin{enumerate}[(ii)]
\item[(ii)] There exist $C_2\in(0,\infty)$ and a nonincreasing function
$h\dvtx[0,\infty)\raw\break  [0,\infty)$ such that $|\phi(\w)| \leq h(\|
\w\|)$ for
all $\|\w\|>C_2$;

\item[(iii)] $\liminf_{\nti} \operatorname{det} (N^{-1}\sum_{k=1}^N G_{\te
_0}(\w
_{kn})G_{\te_0}(\w_{kn})' ) >0$;

\item[(iv)] $\int G_{\te_0}(\w)G_{\te_0}(\w)' \phi^2(\w) \,d\w$ is nonsingular.
\end{enumerate}
\item[(C.4)] (i) $0<\ka< \eta< 1$; and
\begin{enumerate}[(ii)]
\item[(ii)] $\Si_n^{-1/2}\sum_{k=1}^N G_{\theta_0}(\w_{kn})I_n^*(\w_{kn})
\cid N(\0, \bbi_p)$.
\end{enumerate}
\item[(C.5)$'$]
For each $n\geq1$, there exists a function $M_n(\cdot)$ such that
\[
\Biggl\|\sum_{k=1}^N G_{\te_0}(
\w_{kn})G_{\te_0}(\w_{kn})'\exp \bigl(\ci
\mathbf{t}'\w_{kn} \bigr) \Biggr\| \leq M_n(
\mathbf{t}) \qquad\mbox{for all } \mathbf{t}\in\bbr^d
\]
and with $d\nu(\mathbf{t},{\mathbf x}) = \|\mathbf{t}\|\gamma_1(\|
\mathbf{t}\|
)^{
{\de}/{(4+\de)
}}f({\mathbf x}) \,d\mathbf{t}\,d{\mathbf x}$ and $\de\in(0,1]$ of
(C.1),
\begin{eqnarray*}
&&\int\int M_n \bigl(\mathbf{t}+a_1[\mathbf{s}+2
\la_n a_2{\mathbf x}+2\la_n a_3
\mathbf{y}] \bigr) \,d\nu( \mathbf{t},{\mathbf x})\,d\nu(\mathbf{s},\mathbf{y})
\\
&&\qquad= o \bigl( b_n^2 c_n^{1-a_1}
\la_n^{1+a_1} \bigr)\qquad\mbox{for all $a_1,a_2,a_3
\in\{0,1\}$.}
\end{eqnarray*}
\end{enumerate}

We comment on the conditions. Conditions (C.0)--(C.1) are standard
moment and mixing conditions on the spatial process $Z(\cdot)$ (cf.
Lahiri \cite{lahiri2003a}), which entail that $Z(\cdot)$ must be
weakly dependent and are used to ensure finiteness of the variance of
the periodogram values [which are themselves quadratic functions of
$Z(\mathbf{s})$], among other things.
See Doukhan \cite{doukhan1994} for process examples fulfilling such
conditions, including Gaussian, linear and Markov random fields.
Also, note that the function $\g_2(\cdot)$ in (C.0) is allowed to grow
to infinity to ensure validity of the results for bonafide strongly
mixing random fields in $d\geq2$ (cf. Bradley~\cite
{bradley1989,bradley1993}). Condition (C.2) specifies the requirements
on the spatial design density $f$. Part~(i) of (C.2) is a smoothness
condition on $f$ while part (ii) requires the characteristic functions
corresponding to the probability densities $f(\cdot)$ and $f^2(\cdot
)/\int f^2$ to decay at the
rate $O(\|\bolds{\omega}\|^{-a_0})$ as $\|\bolds{\omega}\|\raw
\infty$. Condition (C.2)
is satisfied [with $a_0=d$ in (ii)] when $f(\cdot)$ is the uniform
distribution on a rectangle of the form $(-s_1,t_1)\times\cdots\times
(-s_d,t_d)$ for some $0<s_i, t_i<1/2$ for all $i=1,\ldots,d$. However,
there exist many
nonuniform densities that also satisfy (C.2) with $a_0=d$.

Condition (C.3) specifies the regularity conditions on the spectral
estimating function $G_{\te_0}$. In addition to the spectral moment
condition \eqref{spec-mc}, parts~(i) and~(ii) of (C.3) provide
sufficient conditions that make the errors of Riemann sum
approximations to the variance integral $\int G_{\te_0}(\w) G_{\te
_0}(\w
)'\times\break \phi(\w)\,d\w$ asymptotically negligible. Conditions (C.3)(iii)
and (iv) provide alternative forms of a sufficient condition that
guarantees nonsingularity of the $p\times p$ matrix $\Si_n$
through a subsequence under PID and for the full sequence
under (a subcase of) the MID asymptotic structure, respectively.
Without these, the degrees of freedom of the limiting
chi-squared distribution of the scaled log-SFDEL ratio statistic
can be smaller than~$p$. It is easy to verify that the examples
presented in Section~\ref{fdel} satisfy condition~(C.3), under
mild conditions on the $\mathbf{h}_i$'s in Example~\ref{ex1}, on the
$\mathbf{t}_i$'s in Example~\ref{ex2}, and on the $\mathbf{h}_i$'s
and the
parametric variogram model $2\gamma(\cdot;\te)$ in Example~\ref{ex3}.

Next, consider condition (C.4). The first part of (C.4)
states the requirements on the SFDEL tuning parameters $\ka$ and
$\eta$ that must be chosen by the user in practice. Note that
$\ka$ and $\eta$ determine a Riemann-sum approximation to
the spectral moment condition \eqref{spec-mc} over the
discrete grid \eqref{grid} where $\kappa$ determines the
grid spacing while $\eta$ determines the range of
the approximating set $[ -C^*\la_n^{\eta-\kappa},
C^*\la_n^{\eta-\kappa}]^d$. Thus, one must choose
these parameters so that the grid spacing is small
and the integral of $G_{\theta_0}\phi$ outside
$[ -C^*\la_n^{\eta-\kappa},
C^*\la_n^{\eta-\kappa}]^d$
is small. On the other end, $\kappa$ needs to satisfy the requirement
$0<\ka<1$ to ensure that the neighboring
frequencies in $\N$ are ``asymptotically distant.''
Section~\ref{sim} gives some specific examples of
the choices of $\kappa$ and $\eta$ in finite
sample applications. As for condition (C.4)(ii),
note that the ``asymptotically distant'' property
of the frequencies in $\N$ renders the summands
in $\sum_{k=1}^N G_{\theta_0}(\w_{kn})I^*_n(\w_{kn})$ approximately
independent, and hence, under suitable normalization, the sum must
have a Gaussian limit. One set of sufficient conditions for the
weak convergence of $\Si^{-1/2}_n\sum_{k=1}^NG_{\te_0}(\w_{kn})
[I_n^* (\w_{kn})- E I_n^* (\w_{kn})]$ to a Gaussian limit
is given by a CLT result in Bandyopadhyay, Lahiri and
Nordman \cite{bandy2013c}
(hereafter referred to as [BLN]).
Alternative sufficient conditions for the CLT in (C.4)(ii)
can also be derived requiring that $Z(\cdot)$ is
a $d$-dimensional linear process, but
we do not make any such structural assumptions on $Z(\cdot)$
here.

Finally, consider condition (C.5)$'$ which will be used \textit{only}
in the MID case (cf. Theorems \ref{thmm-2} and \ref{thmm-4}).
This condition can be verified easily when the Fourier transform
$\xi_{j,k}$ (say) of the function $G_{j,\theta_0}G_{k,\theta_0}$ decays
quickly, for all $1\leq j,k\leq p$. In contrast, if
the functions $\xi_{j,k}$ do not decay fast enough, one can
verify (C.5)$'$ using Lemma~\ref{L73} and the arguments in the
proof of the result below, which shows
that condition (C.5)$'$ holds for Examples \ref{ex1}--\ref{ex3}.

%
%
\begin{prop}
\label{prop-1}
For $G_\te(\cdot)$ of Examples \ref{ex1}--\ref{ex3},
condition \textup{(C.5)}$'$ holds.
\end{prop}

The next section states the main results of the paper under PID and MID.
%
\section{Asymptotic distribution of the spatial FDEL ratio statistic}
\label{asy}
\subsection{Results under the PID asymptotic structure}
\label{pid-results}
Let $\px$ denote the joint distribution of the random vectors $\mathbf{X}
_1,\mathbf{X}_2,\ldots,$ generating the locations of the data sites (cf.
Section~\ref{design}). The following result gives the asymptotic
distribution of the SFDEL ratio statistic under PID.

%
%
\begin{thmm}
\label{thmm-1}
Suppose that conditions \textup{(C.0)--(C.4)}
and that $n/\la_n^d \raw c_*\in(0,\infty)$. Then
$
- {}\log\R_n(\te_0) \stackrel{d}{\longrightarrow}\chi^{2}_{p},
\mbox{a.s. ($P_\mathbf{X}$)}$.
\end{thmm}

Theorem~\ref{thmm-1} shows that under conditions (C.0)--(C.4)
[and without requiring~(C.5)$'$], the SFDEL
log-likelihood ratio statistic has an asymptotic chi-squared
distribution, for almost all realizations of the sampling design
vectors $\{\mathbf{X}_i\}$. Note that the scaling for the log
SFDEL ratio statistic is nonstandard---namely, the chi-squared
limit distribution is attained by\break $- {}\log\R_n(\te_0)$, but
\textit{not} by the more familiar form $- {2}\log\R_n(\te_0)$ as
in Wilks' theorem and as in the time series FDEL case
(cf. Nordman and Lahiri \cite{nordman2006}). This is a consequence of
the nonstandard behavior of the periodogram for \textit{irregularly}
spaced spatial data
%
(cf. Section~\ref{fdel}).
However, as
the limit distribution of the SFDEL ratio statistic does not
depend on any unknown population quantities, it
can be used to construct valid large sample tests and
confidence regions for the spectral parameter $\theta$.
Specifically,
a valid large sample level $\al\in(0, 1/2)$ SFDEL test
for testing
%
%
\begin{equation}
H_0\dvtx\te= \te_0\qquad\mbox{vs}\qquad H_1
\dvtx\te\neq\te_0 \label{t-hyp}
\end{equation}
will reject $H_0$ if $ - {} \log\R_n(\te_0) > \chi^2_{1-\al,p}$,
where $\chi^2_{1-\al,p}$ denotes the $(1-\al)$ quantile
of the $\chi^2_{p}$-distribution. For SFDEL based confidence regions
for $\te$,
a similar distribution-free calibration holds (cf. Section~\ref{unif}).

\begin{remark}
Note that the distribution of
$\R_n(\te_0)$ depends on two sources of randomness, namely,
the spatial process $\{Z(\mathbf{s}) \dvtx\mathbf{s}\in{\mathbb
{R}^{d}}\}$ and the
vectors $\{\mathbf{X}_i\}_{i \geq1}$. Let $\L(T|{\cal X})$ denote
the conditional distribution of a random variable (based on
both $\{Z(\cdot)\}$ and $\{\mathbf{X}_i\}$), given ${\cal X}\equiv
{\sigma}\langle\mathbf{X}_1, \mathbf{X}_2,\ldots\rangle$ and let $d_L$
denote the Levy metric on the set of probability distributions
on $\bbr$.
Then
a more precise statement of the Theorem~\ref{thmm-1} result, under the
conditions given there, is
\[
d_L \bigl( \L \bigl( -\log\R_n(\theta_0) |{
\cal X} \bigr), \chi^2_p \bigr) =o(1) \qquad\mbox{a.s. ($
\px$).}
\]
A similar interpretation applies to the other theorems
presented in the paper.
\end{remark}

\subsection{Results under the MID asymptotic structure}
\label{mid-results}
The limit behavior of the spatial FDEL ratio statistic under
the MID asymptotic structure shows a more complex pattern and it
depends on the strength of the infill
component. Note that $c_n = n/\la_n^d$ denotes the
relative growth rate of the sample size and the volume of the sampling
region of $\D_n$, and hence,
$
c_n\raw\infty$ as $\nti$
{under MID},
with a higher the value of $c_n$ indicating a
higher rate of infilling. The following result gives the
asymptotic behavior of the SFDEL ratio statistic under different
growth rates of~$c_n$.

%
%
\begin{thmm}
\label{thmm-2}
Suppose that conditions \textup{(C.0)--(C.4)} and \textup{(C.5)}$'$
hold
[where \textup{(C.3)} may be replaced by \textup{(C.3)(i), (ii), (iv)}
for part \textup{(b)}].
\begin{itemize}[(a)]
\item[(a)] \textsc{(MID with a slow rate of infilling)}.
If $1\ll c_n^2 \ll N\la_n^{-\ka d}$, then
$
- {}\log\R_n(\te_0) \stackrel{d}{\longrightarrow}\chi^{2}_{p},
\mbox{a.s. ($P_\mathbf{X}$).}
$
\item[(b)]
\textsc{(MID with a fast rate of infilling)}.
If $c_n^2 \gg N\la_n^{-\ka d}$, then\break
$- {2}\log\R_n(\te_0) \stackrel{d}{\longrightarrow}\chi^{2}_{p},
\mbox{a.s. ($P_\mathbf{X}$).}
$
\end{itemize}
\end{thmm}

Theorem~\ref{thmm-2} shows that the asymptotic distribution
of $-\log\R_n(\te_0)$ can be different depending on the rate
at which the infilling factor $c_n$ goes to infinity. When the
rate of decay in $c_n^2$ is slower than the critical rate
$N\la_n^{-\ka d}\sim\la_n^{(\eta-\ka)\,d}$,
corresponding to the asymptotic volume of the frequency grid (i.e.,
determined by the number $N\propto\la_{n}^{d\eta}$ of frequency points
on a regular grid and the ${\mathbb{R}^{d}}$-volume $\lambda
_n^{-d\kappa}$ between
grid points), the negative log SFDEL ratio
has the same limit distribution as in the PID case. However,
when the factor $c_n^2$ grows at a faster rate than
$ \la_n^{(\eta-\ka)d}$, the \textit{more familiar} version
of scaling $-2$ is appropriate
for the log SFDEL ratio. From the proof of Theorem~\ref{thmm-2},
it also follows that
in the boundary case, that is, when
$c_n^2 \sim\la_n^{(\eta-\ka) d}$,
the limit distribution of $-\log\R_n(\te_0) $ is determined
by that of a quadratic form in independent Gaussian random
variables and is \textit{not} distribution-free. As a result, this
case is not of much interest from an applications
point of view. However, as $c_n = n/\la_n^d$ is a \textit{known}
factor, one can always choose the SFDEL tuning parameters
$\ka, \eta$ to avoid the boundary case. 

\begin{remark}
Theorem~\ref{thmm-2} shows that when
the rate of infilling $c_n$ does not grow too fast,
the presence of the infill component does not have
an impact on the asymptotic distribution of the log SFDEL
ratio statistic. Thus, the limit behavior under the PID asymptotic
structure has a sort of robustness that extends beyond its realm
and covers parts of the MID asymptotic structure in the frequency
domain. This is very much different from the known results on the
limit distributions of the sample mean and of asymptotically
linear statistics in the spatial domain where \textit{all subcases} of
the MID asymptotic structure lead to the \textit{same} limit
distribution and where the MID limit is \textit{different}
from the limit distribution in the PID case
(cf. Lahiri \cite{lahiri2003a}, Lahiri and Mukherjee \cite{lahiri2004}).
\end{remark}

\subsection{A unified scaled spatial FDEL method}
\label{unif}
Results of Sections \ref{pid-results} and \ref{mid-results} show that
in the spatial case, the standard calibration of the EL ratio statistic
may be incorrect depending on the relative rate of infilling. Although
nonstandard, $- {2}\log\R_n(\te_0)$ has the same ${2}\chi^{2}_{p}$
distribution under the PID spatial asymptotic structure for all values
of $c_*=\lim_{\nti}n/\la_n^d$. In contrast, the limit distribution
of $
- {2}\log\R_n(\te_0)$ can change from the nonstandard ${2}\chi
^{2}_{p}$ to the standard $\chi^{2}_{p}$ under the MID asymptotic
structure when the rate of infilling is faster.
While this gives rise to a clear dichotomy in the limit, the choice of
the correct scaling constant, and hence, the correct calibration may
not be obvious in a finite sample application. To deal with this
problem, we develop a data based scaling factor that adjusts itself to
the relative rates of infilling and delivers a unified $\chi^{2}_{p}$
limit law under the PID as well as under the different subcases of the
MID. Specifically, define the modified FDEL statistic
\[
- 2 a_n(\te) \log\R_n(\te),
\]
where
%
%
\begin{equation}
a_n (\te) = \frac{\sum_{k=1}^N \|G_{\te}(\w_{kn})\|^2 \tI_n^2(\w_{kn})}{
\sum_{k=1}^N \|G_{\te}(\w_{kn})\|^2 {I}_n^2(\w_{kn})
}. \label{an1}
\end{equation}
Note that for any $\te$, the factor $a_n(\te)$ can be computed using
the data $\{Z(\mathbf{s}_1),\ldots,  Z(\mathbf{s}_n)\}$, where the
numerator of
$a_n(\theta)$ is computed using the bias-corrected periodogram while the
denominator is based on the raw periodogram. For the testing problem
$H_0\dvtx\te=\te_0$ against $H_1\dvtx\te\neq\te_0$, this
requires computing
the factor $a_n(\cdot)$ once. However, for constructing confidence
intervals, $a_n(\te)$ must be computed repeatedly and, therefore, this
version of the SFDEL is somewhat more computationally intensive.

To gain some insight into the choice of $a_n(\te)$, note that it is
based on the ratio of the sums of the periodogram and its
bias-corrected version that are weighted by the squared norms of the
function $G_{\te}(\cdot)$ at the respective frequencies $\w_{kn}$. As
explained before, the bias correction of the periodogram of irregularly
spaced spatial data is needed to render the EL-moment condition in
\eqref{Fdel} unbiased. However, this leads to a ``mismatch'' between the
variance of the sum $\sum_{k=1}^N G_{\te}(\w_{kn}) \tilde{I}_n(\w
_{kn})$ and the automatic scale adjustment factor provided by the EL
method. The numerator and the denominator of $a_n(\theta)$ capture the
effects of this mismatch under different rates of infilling and hence,
$a_n(\cdot)$ provides the ``correct''
scaling constant under the different asymptotic regimes considered here.

We have the following result on the modified SFDEL ratio statistic.

%
%
\begin{thmm}
\label{thmm-4} Suppose that the conditions of one of Theorems \ref
{thmm-1}--\ref{thmm-2}
hold. Then, under $\te=\te_0$,
%
%
\begin{equation}
- 2 a_n(\te_0) \log\R_n(\te_0)
\stackrel{d} {\longrightarrow}\chi^{2}_{p}\qquad\mbox{a.s. ($P_\mathbf{X}$).}
\end{equation}
\end{thmm}

Theorem~\ref{thmm-4} shows that the modified SFDEL method can
be calibrated using the quantiles of the chi-squared
distribution with $p$ degrees of freedom \textit{for all}
of the three asymptotic regimes covered by Theorems
\ref{thmm-1}--\ref{thmm-2}.
Thus, the empirically scaled log-SFDEL ratio
statistic provides a unified way of testing and constructing
confidence sets under different asymptotic regimes.
Specifically, for any $\al\in(0,1/2)$,
\[
\C_\al\equiv \bigl\{\te\in\Th\dvtx- {2a_n(\theta)} \log
\R_n(\te) \leq\chi^2_{1-\al,p} \bigr\}
\]
gives a confidence region for the unknown parameter $\te$
that attains the nominal confidence level $(1-\al)$
asymptotically. 
The main advantage of the SFDEL method here is that
we do \textit{not} need to find a studentizing covariance matrix
estimator explicitly, which by itself
is a nontrivial task, as this would require \textit{explicit estimation}
of the spectral density $\phi(\cdot)$ and the
spatial sampling density $f(\cdot)$ under different asymptotic regimes.

\section{Numerical results}
\label{sim}
\subsection{Results from a simulation study}\label{sec6.1}
%
Here, we examine the coverage accuracy of the SFDEL method
in finite samples,
applied to a problem of variogram model fitting described in
Section~\ref{examples}. We consider an exponential variogram
model form (up to variance normalization)
\[
2\gamma(\mathbf{h}; \theta_1,\theta_2) = 1 - \exp \bigl[-
\theta_1 |h_1| - \theta_2 |h_2|
\bigr],
\]
with parameters $\theta_1,\theta_2>0$ where $\mathbf
{h}=(h_1,h_2)^\prime
\in\mathbb{R}^2$. Over several sampling region sizes $\D_{n}
=\lambda_n[-1/2,1/2)^2 $, $\lambda_n=12,24,48$, and sample sizes
$n=100,\break 400,900,1400$, we generated i.i.d. sampling sites
$\mathbf{s}_1,\ldots,\mathbf{s}_n\in\D_n$ and real-valued
stationary Gaussian responses $Z(\cdot)$ following the exponential
variogram form with $\theta_1=\theta_2=1$ and $E Z(\mathbf{s})=0$,
$\operatorname{Var}[Z(\mathbf{s})]=1$ (the simulation results are
invariant here
to values for the mean and variance). In the spatial sampling
design (cf. Section~\ref{design}), two distributions $f$
for sites were considered, one being uniform over $\D_0$ and
the other being a mixture
of two bivariate normal distributions
$0.5N((0, 0)^\prime, \bbi_2)+0.5N((1/4, 1/4)^\prime, 2\bbi_2)$,
truncated outside $\D_0$,
where $\bbi_2$ denotes a $2\times2$ identity matrix.

In implementing the modified SFDEL method to compute 90\% confidence
regions for $\theta= (\theta_1,\theta_2)$, we used the estimating
functions in (\ref{vrg-eeq}) over $m=2$ sets of lags $\mathbf{h}_1,
\mathbf{h}_2 \in\mathbb{R}^2$ and evaluated the (sample
mean-centered) periodogram at scaled frequencies
$ \mathcal{N}_n=\{ \lambda_n^{-\kappa} \mathbf{j} \dvtx\mathbf
{j}\in\mathbb{Z}^2 \cap[-C^* \lambda_n, C^* \lambda_n]^2\}$;
we varied values $C^*=1,2,4$ (with $\eta=1$ held fixed)
and $\kappa=0.05,0.1,0.2$ along with
considering different combinations of lags $\mathbf{h}_1, \mathbf{h}_2$.
Recall that $C^*$ and $\kappa$ respectively control the number
and spacing of periodogram ordinates, where choices of $\kappa$
here roughly induce spacings between frequencies of 1, 0.75 or 0.5
in horizontal/vertical directions; in our findings, these spacings
were adequate whereas tighter spacings (e.g., $\kappa\geq0.4$)
tended to perform less well by inducing stronger dependence
between periodogram ordinates.

%
%
\begin{table}
\tabcolsep=0pt
\caption{Coverage percentage of 90\% SFDEL regions for variogram model
parameters $\theta$ (uniform design)}
\label{vario-results}
%
\begin{tabular*}{\textwidth}{@{\extracolsep{\fill}}ld{1.2}cccccccccccc@{}}
\hline
&&\multicolumn{12}{c@{}}{$\mathbf{{h}_1}\bolds{=(1,1)^\prime,}
\mathbf{{h}_2}\bolds{=(1,-1)^\prime}$}\\[-6pt]
&&\multicolumn{12}{c@{}}{\hrulefill}\\
&& \multicolumn{4}{c}{$\bolds{\lambda_n=12}$}&
\multicolumn{4}{c}{$\bolds{\lambda_n=24}$}&
\multicolumn{4}{c@{}}{$\bolds{\lambda_n=48}$}\\[-6pt]
&& \multicolumn{4}{c}{\hrulefill}&
\multicolumn{4}{c}{\hrulefill}&
\multicolumn{4}{c@{}}{\hrulefill}\\
\multicolumn{1}{@{}l}{$\bolds{C^{*}}$}&
\multicolumn{1}{c}{$\bolds{\kappa}$} &
\multicolumn{1}{c}{$\mathbf{100}$} & \multicolumn{1}{c}{$\mathbf
{400}$} &
\multicolumn{1}{c}{$\mathbf{900}$} & \multicolumn{1}{c}{$\mathbf
{1400}$} &
\multicolumn{1}{c}{$\mathbf{100}$} & \multicolumn{1}{c}{$\mathbf
{400}$} &
\multicolumn{1}{c}{$\mathbf{900}$} & \multicolumn{1}{c}{$\mathbf
{1400}$} &
\multicolumn{1}{c}{$\mathbf{100}$}
& \multicolumn{1}{c}{$\mathbf{400}$} & \multicolumn{1}{c}{$\mathbf
{900}$} &
\multicolumn{1}{c@{}}{$\mathbf{1400}$}\\
\hline
1&0.05&86.4&85.6&82.0&80.3&88.9&87.8&87.8&89.9&89.3&89.4&89.7&87.9\\
1&0.1&87.1&85.3&78.6&75.9&89.0&90.2&89.6&90.4&89.0&91.4&91.5&90.0 \\
1&0.2&86.5&85.1&81.1&76.4&90.0&88.7&90.1&89.7&87.6&87.9&87.9&88.9 \\[3pt]
2&0.05&88.1&87.8&86.1&85.9&89.0&88.6&89.7&87.9&89.2&88.9&90.5&89.7 \\
2&0.1&86.6&86.8&86.2&84.2&89.2&88.4&91.1&89.9&90.6&90.0&90.0&91.4 \\
2&0.2&89.6&88.8&84.6&83.8&88.9&89.9&89.9&89.2&89.9&89.3&88.1&89.4 \\[3pt]
4&0.05&89.0&87.8&89.6&88.1&89.3&89.0&90.1&90.2&92.9&88.2&90.6&89.9 \\
4&0.1&86.3&88.6&88.7&86.4&90.3&89.4&90.3&89.2&92.0&87.8&90.8&89.1 \\
4&0.2&88.4&89.0&87.4&87.9&88.7&88.9&90.0&89.6&92.8&88.6&88.5&88.8 \\
\hline
\end{tabular*}
%
\end{table}

%
%
\begin{table}[b]
\tabcolsep=0pt
\caption{Coverage percentage of 90\% SFDEL regions for variogram model
parameters $\theta$ (nonuniform~design)}
\label{vario-results2}
\begin{tabular*}{\textwidth}{@{\extracolsep{\fill}}lccccccccccccc@{}}
\hline
&&\multicolumn{12}{c@{}}{$\mathbf{{h}_1}\bolds{=(1,1)^\prime,}
\mathbf{{h}_2}\bolds{=(1,-1)^\prime}$}\\[-6pt]
&&\multicolumn{12}{c@{}}{\hrulefill}\\
&& \multicolumn{4}{c}{$\bolds{\lambda_n=12}$}& \multicolumn
{4}{c}{$\bolds{\lambda
_n=24}$}& \multicolumn{4}{c@{}}{$\bolds{\lambda_n=48}$}\\[-6pt]
&& \multicolumn{4}{c}{\hrulefill}&
\multicolumn{4}{c}{\hrulefill}&
\multicolumn{4}{c@{}}{\hrulefill}\\
\multicolumn{1}{@{}l}{$\bolds{C^{*}}$}&
\multicolumn{1}{c}{$\bolds{\kappa}$} &
$\mathbf{100}$ & $\mathbf{400}$ & $\mathbf{900}$ & $\mathbf{1400}$ &
$\mathbf{100}$ & $\mathbf{400}$ & $\mathbf{900}$ & $\mathbf{1400}$
& $\mathbf{100}$
& $\mathbf{400}$ & $\mathbf{900}$ & \multicolumn{1}{c@{}}{$\mathbf
{1400}$}\\
\hline
1&0.05&88.3&86.8&85.5&79.6&89.4&88.9&86.4&90.1&90.2&89.2&89.6&90.0\\
1&0.10&85.7&83.5&80.3&78.7&88.8&87.7&89.6&92.0&87.0&90.5&89.5&89.3\\
1&0.20&87.9&86.0&82.4&79.1&89.6&90.0&90.7&90.0&87.4&88.2&88.8&88.7 \\[3pt]
2&0.05&89.4&89.3&88.0&83.8&90.1&88.6&89.7&88.9&89.6&90.7&89.5&91.2 \\
2&0.10&86.2&87.7&84.3&85.9&89.0&90.7&90.0&88.4&90.1&91.5&90.1&90.0 \\
2&0.20&88.7&89.5&88.5&85.6&90.7&90.4&89.7&88.5&89.7&88.5&90.3&89.8 \\[3pt]
4&0.05&89.5&89.5&88.3&88.0&87.7&90.0&88.6&90.7&91.8&89.8&89.2&90.9 \\
4&0.10&87.0&88.8&87.4&86.2&89.0&89.9&87.9&89.4&91.7&87.9&88.2&89.9 \\
4&0.20&90.6&89.1&89.1&86.3&89.5&89.0&87.9&89.4&91.5&90.2&89.3&90.1 \\
\hline
\end{tabular*}
\end{table}

The coverage results (based on
1000 simulation runs) are listed in Tables \ref{vario-results}--\ref
{vario-results2} for the lag
$\mathbf{h}_1=(1,1)^\prime, \mathbf{h}_2=(1,-1)^\prime$
for the uniform and
nonuniform spatial sites, respectively, with the
results for the other sets of lags reported in
the supplementary material \cite{supp}.
Except for the occasions with the smallest lag combination
[$\mathbf{h}_1=(1,1)^\prime, \mathbf{h}_2=(1,-1)^\prime$] and
the smallest sampling region $\lambda_n=12$ with large $n$, the
coverages tended to agree quite well with the nominal level.
Further, the coverage levels were largely insensitive to the
number and spacing of periodogram ordinates for various
sample and region sizes. Results for both stochastic sampling
designs were also qualitatively similar.
%

\subsection{An illustrative data example}\label{sec6.2}

%
%
\begin{figure}

\includegraphics{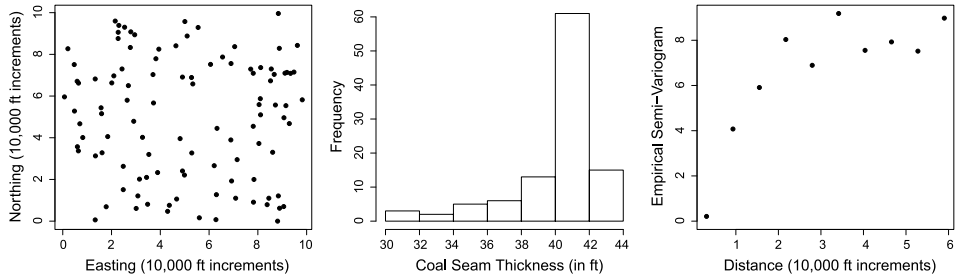}

\caption{Coal seam data: Sampling locations, distribution of thickness,
empirical semivariogram.}\label{fig1}
\end{figure}

As a brief demonstration of the SFDEL method, we consider a coal seam
dataset based on a SAS example (\cite{SAS2008}, Chapter~70).
Figure~\ref{fig1} shows locations of 105 sampling sites and the corresponding
distribution of coal seam thickness.
Coal seam measurements often exhibit spatial smoothness
(Journel and Huijbregts \cite{journel1978}, page 165), as also indicated
in the empirical semivariogram in Figure~\ref{fig1}
(found by binning distances into 10 bins up to half the maximum
distance between points
and plotting Matheron's average over each bin against the bin
midpoint). Following the SAS analysis, this suggests
a Gaussian variogram model $2\gamma(\mathbf{h}; \theta_1,\theta_2) =
2\theta_1 [1 - \exp( -\|\mathbf{h}\|^2/\theta_2^2)]$, $\mathbf
{h}\in
\mathbb{R}^2$,
with scale $\theta_1>0$ and range $\theta_2>0$ parameters, though the present
data are synthetic with a value $\theta_2=1$ as explained below.

Note that the spatial locations are not clearly uniform nor is the
marginal distribution apparently normal.
To fit the variogram model in a way that allows nonparametric
confidence intervals (CIs) to access
the precision of the estimated parameters, without making assumptions
about the joint distribution of the data
or the distribution of spatial locations, one can apply the SFDEL
method using estimating functions as in Example~\ref{ex3} motivated
by least squares estimation.
Alternatively, one can apply a kernel bandwidth estimator of the
varigoram for which
large sample distributional results
are recently known (Garc\'{i}a-Soid{\'a}n \cite{gsoidan2007}, Maity and
Sherman~\cite{maity2012}).

We focus on the range parameter $\theta_2$.
Using a lag set $\mathbf{h}_1=(1/4,1/4)^\prime$, $\mathbf
{h}_2=(1,1)^\prime$,
$\mathbf{h}_3=(2,2)^\prime$ in SFDEL, motivated by empirical lags in
Figure~\ref{fig1}, the maximized
SFDEL function produces a point estimate $\hat{\theta}_2 = 1.123$
($\times10\mbox{,}000$ ft)
with a 90\% SFDEL CI for $\theta_2$ as $(0.896, 1.571)$.
This arises from a frequency grid
$\{\lambda_n^{-\kappa} \mathbf{j} \dvtx[ -C^{*}\lambda_n,
C^{*}\lambda
_n]^2 \cap\mathbb{Z}^2\}$, $C^{*}=2$, $\kappa=0.2$
based on $\lambda_n=10$ for sampling region in
Figure~\ref{fig1}. With a larger frequency grid $C^{*}=4$, $\kappa=0.2$,
the 90\% SFDEL CI is similar $(0.887, 1.378)$ with a point estimate
$1.071$, and increasing the
grid spacing $\kappa=0.1$ produces similar range estimates ($1.107$
and $1.101$ for $C^{*}=2,4$) and
intervals. In contrast, using the lags above and the Nadaraya--Watson
kernel estimator of the semivariogram
$\hat{\gamma}(\mathbf{h}) $
(based on the Epanechnikov kernel, cf. Garc\'{i}a-Soid{\'a}n \cite{gsoidan2007}),
the range parameter estimates are $ 1.505, 1.230, 1.335$ for
bandwidths $h=0.5,1,1.5$, where $h=0.5$ arises from MSE optimal order
considerations.
This approach can also produce large-sample nonparametric CIs based on
normal limits for
$\lambda_n[\hat{\gamma}(\mathbf{h})- \gamma(\mathbf{h})]$, having a
covariance matrix $\mathcal{C}\cdot V$, $\mathcal{C}=[\int f^2]^{-2}
\int f^4$, involving
the unknown density $f$ of locations $\{\mathbf{s}_{i}/\lambda_n\}
_{i=1}^{105}$ on $[0,1]^2$.
For bandwidths $h=0.5,1,1.5$, the 90\% CIs for $\theta_2$ are given by
$(1.203, 1.743)$, $(0.972, 1.512)$, $(1.038, 1.564)$
based on $\hat{\mathcal{C}}=1.23$ from bivariate kernel density
estimation (cf. Venables and Ripley \cite{ripley2002}) and simplifying
the matrix $V$ by assuming
the process is Gaussian (cf. Garc\'{i}a-Soid{\'a}n \cite{gsoidan2007},
pages~490--491). Unlike CIs from kernel estimation, the SFDEL CIs
require no variance or density approximation steps, tend to be less
sensitive to tuning parameters,
and all contain the true value $\theta_2=1$ here.

To provide some assessment of the CI methods, we conducted a small
simulation study
generating marginally standard normal variates $\{Y(\mathbf{s}_i)\}
_{i=1}^{105}$
with correlation $\operatorname{corr}[Y(\mathbf{s}),Y(\mathbf
{s}+\mathbf
{h})]=\exp(-\|\mathbf{h}\|^{2})$ at the locations $\{\mathbf{s}_i\}
_{i=1}^{105}$ in Figure~\ref{fig1} and
defining observations $\{Z(\mathbf{s}_i)= \sqrt{\theta
_1/2}[Y^2(\mathbf
{s}_i)-1] + 40.23 \}_{i=1}^{105}$ from a spatial process
having a Gaussian variogram as above with scale $\theta_1=7.5$
and range $\theta_2=1$; this data-generation approximately
matches features in the original SAS coal seam data and also produced
the data example above.
Based on 1000 simulations, 90\% CIs for the range parameter $\theta_2$
from the SFDEL method had coverages $90.5, 87.4, 93.3$ for $C^{*}=2$
and $89.3, 88.7, 86.4$ for $C^{*}=4$, over grid
spacings $\kappa=0.05,0.1,0.2$. In contrast,
90\% CIs for $\theta_2$ from the kernel estimation approach had actual
coverages $68.4, 74.4, 52.7$ for bandwidths $h=0.5,1,1.5$.

\section{Proofs of the results}
\label{proof_fdel}
%
\subsection{Notation and lemmas}\label{sec7.1}
Define the bias corrected periodogram\break $\tilde{I}_{n}(\w)={I}_{n}(\w) -
n^{-1}\lambda_n^d \hsin(\0)$ and its (unobservable) variant
${I}^*_{n}(\w)={I}_{n}(\w) - n^{-1}\lambda_n^d {\sigma}(\0)$.
Recall that $A_n(\bolds{\omega}) = c_n^{-1} {\sigma}(\0) +
K \phi(\bolds{\omega})$, $\bolds{\omega}\in\bbr^d$, where
$K=(2\pi)^d \int f^2$.
For notational simplicity, for a random quantity $T$ depending on both
$\{Z(\mathbf{s})\dvtx \mathbf{s}\in\bbr^d\}$ and $\{\mathbf
{X}_1,\mathbf{X}_2,\ldots\}$,
$ET$ will
denote the conditional expectation of $T$ given
$\bbX\equiv\{\mathbf{X}_1,\mathbf{X}_2,\ldots\}$ and likewise $P$
will denote
conditional probability. Thus, in the following,
$P(-2\log\R_n(\theta_0) \leq t)$ in fact refers to
\[
P \bigl(-2\log\R_n(\theta_0) \leq t | \bbX \bigr),\qquad
t>0.
\]
Also, write
$P_{\mathbf{X}}$ and $E_{\mathbf{X}}$ to denote the probability and
the expectation
under the joint distribution of $\mathbf{X}_1,\mathbf{X}_2,\ldots.$
Further, let
$C$ or $C(\cdot)$ denote generic constants that depend on their
arguments (if any), but do not depend on $n$ or the $\{\mathbf{X}_i\}$.

We now provide some lemmas that will be used for proving
the main results of the paper. {Proofs of the
lemmas and Proposition~\ref{prop-1} are relegated to the
supplementary material \cite{supp}}
to save space.
For continuity, supplementary material \cite{supp}\vadjust{\goodbreak} begins with three technical lemmas
(Lemmas 7.1--7.3), providing
some general
cumulant and integral inequalities as well as the bias and the variance
of the
periodogram $I_n(\cdot)$ that are used to establish Lemmas \ref
{L73}--\ref{L77}
below; these
results may also be of independent
interest. As presented next, Lemmas \ref{L73}--\ref{L77} deal
with various
properties and sums of the periodogram
that we will need to analyze the asymptotic behavior of the SFDEL ratio
statistic under different asymptotic structures and establish the main
results in Section~\ref{sec7.2}.

\newtheorem{lemma}[thmm]{Lemma}

\setcounter{thmm}{3}
\begin{lemma}
\label{L73}
Under conditions \textup{(C.0)--(C.3)} and \textup{(C.5)}$'$
\begin{eqnarray*}
&&E \Biggl[ \sum_{k=1}^N
G_{\te_0}(\w_{kn}) G_{\te_0}(\w_{kn})'
I_n^2(\w_{kn}) \Biggr]
\\
&&\qquad= 2 \sum_{k=1}^N
G_{\te_0}( \w_{kn}) G_{\te_0}(\w_{kn})'
A_n^2(\w_{kn}) + o \bigl(b_n^2
\bigr) \qquad\mbox{a.s. ($\px$).}
\end{eqnarray*}
\end{lemma}

%
\begin{lemma}
\label{L75}
Under conditions \textup{(C.0)--(C.3)} and \textup{(C.5)}$'$,
\begin{eqnarray*}
&&\sum_{i=1}^{N}G_{\theta_{0}}( \bolds{
\omega}_{in})G_{\theta_{0}}(\bolds{\omega}_{in})^{\prime}
\bigl[ \tI_{n}^{2}(\bolds{\omega}_{in}) -
\bigl(A_n^2(\w_{in}) + K^2
\phi^2(\w_{in}) \bigr) \bigr]
\nonumber
\\[-8pt]
\\[-8pt]
\nonumber
&&\qquad= o_p
\bigl(b_n^2 \bigr)\qquad%
\mbox{a.s. ($\px$).}
\end{eqnarray*}
\end{lemma}

%


%
\begin{lemma}
\label{L76}
Under conditions \textup{(C.0)--(C.3)} and \textup{(C.5)$'$},
for any $\varepsilon>0$,
$
P (\max_{1\leq k\leq N} \llVert G_{\theta_{0}}(\bolds{\omega}_{kn})
I_{n}(\bolds{\omega}_{kn})\rrVert>\ep b_n )=o(1), \mbox{a.s. }
(P_{\mathbf{X}})$.
\end{lemma}

\begin{lemma}
\label{L77}
Let $ch^o(B)$ denote the interior of the convex hull of a set $B
\subset\bbr^p$. Under conditions \textup{(C.0)--(C.3)} and \textup
{(C.5)$'$}, it holds
that, as
$n\rightarrow\infty
$,
$
P (0\in\mbox{ch}^{o}\{G_{\theta_{0}}(\bolds{\omega}_{kn})
\tilde{I}_{n}(\bolds{\omega}_{kn})\}_{k=1}^{N} )\rightarrow1$
$\mbox
{a.s. }
(P_{\mathbf{X}})$.
\end{lemma}
%

\subsection{Proofs of the main results}\label{sec7.2} We now present
the proofs of the results from Section~\ref{asy}. In the following, references
to the equations
from the supplementary material \cite{supp} are given as (S.$*$).

\begin{pf*}{Proof of Theorem~\ref{thmm-1}}
By Lemma~\ref{L77}, $\R_{n}(\theta_{0})$ exists
and is positive on a set with probability tending to one,
a.s. ($\px$).
When $\R_{n}(\theta_{0})>0$ holds, by a general and standard EL result
based on Lagrange multipliers
(cf. Owen \cite{owen1988}, page 100), one can express $\R_n(\theta
_0)$ as
%
%
\begin{equation}
\R_{n}(\theta_{0})=\prod_{k=1}^{N}(1+
\g_{k})^{-1}, \label{el-1}
\end{equation}
where $\b_{\theta_{0}}\equiv\b_{\theta_{0}, n}$ satisfies
$F_{n}(\theta_{0},\b_{\theta_{0}})=0$ for
\[
F_{n}(\theta,\b)\equiv N^{-1}\sum
_{k=1}^{N}\frac{G_{\theta}(\w_{kn})
\tI_{n}(\w_{kn})}{1+\b^{\prime}G_{\theta}(\w_{kn})\tI_{n}(\w_{kn})}
\]
and where
$\g_{k}\equiv\g_{k,n}=\b_{\theta_{0}}^{\prime}G_{\theta_{0}}(\w
_{kn}) \tilde
{I}_{n}(\w_{kn})$
satisfies $|\g_{k}|<1$ for all
$1\leq k \leq N$.
To prove the theorem, it is enough
to show that given any subsequence $\{n_i\}$, there exists a
further subsequence $\{n_k\}$ of
$\{n_i\}$ such that $-\log\R_{n_k}(\te_0) \stackrel{d}{\rightarrow}
\chi^2_p$.
We use this line of argument because, as the proof indicates, the asymptotic
expansion of $-\log\R_{n_k}(\te_0)$ involves mean-like quantities
(e.g., term $J_k$ in the following)
which may have differing (normal) limit distributions along different
subsequences of $\{n_k\}$; nevertheless,
the log-ratio statistic $-\log\R_{n}(\te_0)$ is shown to have a
single, well-defined chi-square limit.

Note first that, under the PID structure here, it follows immediately
from (C.3) [cf. (S.10)] that
%
%
\begin{equation}\qquad
b_n^2\sim Nc_*^{-2}\quad\mbox{and}\quad \Biggl\|
\Si_n - 2c_*^{-2} \bigl[{\sigma}(\mathbf{0})
\bigr]^2 \sum_{k=1}^NG_{\te
_0}(
\w_{kn}) G_{\te_0}(\w_{kn})' \Biggr\| = o(N)
\label{pid-1}
\end{equation}
[which is applied to show (\ref{JW-bd}) from (\ref{pid-2}) next].
Fix a subsequence $\{n_i\}$. Then by (C.3)(iii) and the fact that
$\|G_{\te_0}(\w)\| \leq C$ for all $\w\in\bbr^d$, it follows that
there exist a
subsequence $\{n_k\}$ of
$\{n_i\}$ and a nonsingular matrix $\Ga^*$ (possibly depending on $\{
n_k\}$) such that
%
%
\begin{equation}
N^{-1}\sum_{j=1}^N
G_{\theta_{0}}(\w_{jn})G_{\theta_{0}}(\w_{jn})'
\raw\Ga^* \qquad\mbox{through $\{n_k\}$.} \label{pid-2}
\end{equation}
For simplicity, replace the subscript $n_k$ by $k$
and set $N_k\equiv N_{n_k}$, $\w_j \equiv\w_{j, n_k}$,
$\g_j\equiv\g_{j,n_k }$ and
$\b_{\te_0}\equiv\b_{ \te_0, n_k}$. Also, let
$
W_{k}
=
N_k^{-1}\sum_{j=1}^{N_k}G_{\theta_{0}}(\w_{j})G_{\theta_{0}}(\w
_{j})^{\prime}
\tI_{k}^{2}(\w_{j})$ and
$J_{k}= N_k^{-1}\sum_{j=1}^{N_k}G_{\theta_{0}}(\w_{kn})\tI_{k}(\w_{j})$.
Then by \eqref{pid-1}, \eqref{pid-2}, condition (C.3) and
Lemmas \ref{L75}--\ref{L76}, we have, a.s. ($\px$),
%
%
\begin{equation}
|J_{k}|= O_{p} \bigl(N_k^{-1/2}
\bigr)\quad\mbox{and}\quad \bigl\llVert W_{k}-(2N_k)^{-1}
\Si_{k} \bigr\rrVert=o_{p}(1) \label{JW-bd}
\end{equation}
as $k\raw\infty$. Since $N_k^{-1}\Si_k \raw2[{\sigma}(\mathbf{0})]^2
c_*^{-2}\Ga^*$,
$ W_{k}$ is nonsingular
whenever $\|W_k - [{\sigma}(\mathbf{0})]^2 c_*^{-2}\Ga^* \|$ is sufficiently
small.

\begin{claim*}\label{claim}
$\llVert\b_{\theta_{0}}\rrVert
=O_{p}(N_k^{-1/2})$, a.s.
($\px$).
\end{claim*}

\begin{pf}
Write $\b_{\theta_{0}}=t_{0}\mathbf{u}_{0}$ where
$\llVert\mathbf{u}_{0}\rrVert=1$ and $t_{0}=\|\b_{\theta_{0}}\|$. Then
\begin{eqnarray*}
0&=& \bigl\llVert F_{k}(\theta_{0},\b_{\theta_{0}})
\bigr\rrVert \geq\bigl\vert\mathbf{u}_{0}^{\prime}F_{k}(
\theta_{0},t_{\theta
_{0}})\bigr\vert
\\
&=& N_k^{-1}\Biggl\vert\mathbf{u}_{0}^{\prime}
\Biggl(\sum_{j=1}^{N_k}G_{\theta
_{0}}(
\w_j) \tI_{k}(\w_{j})-t_{0}\sum
_{j=1}^{N_k}\frac{ G_{\theta_{0}}(\w_{j})
\tI_{k}(\w_{j})\mathbf{u}_{0}^{\prime}G_{\theta_{0}}(\w_{j})
\tI_{n}(\w_{j})}{1+t_{0}\mathbf{u}_{0}^{\prime} G_{\theta_{0}}(\w_{j})
\tI_{n}(\w_{j})} \Biggr)\Biggr\vert
\\
&\geq& \frac{t_{0}\mathbf{u}_{0}^{\prime}W_{k}\mathbf{u}_{0}}{1+t_{0}Y_{k}} -\sum_{j=1}^{p}
\bigl\vert\mathbf{e}_{j}^{\prime}J_{k}\bigr\vert,
\end{eqnarray*}
where $Y_{k}= \max_{1\leq j\leq N_k}
\llVert G_{\theta_{0}}(\w_{j})\rrVert |\tI_{k}(\w_{j})|$ and
$\mathbf{e}
_1,\ldots,\mathbf{e}_r$ denote the standard basis of $\bbr^r$, with
$\mathbf{e}
_i\in\bbr^r$ having a $1$ in the $i$th position and $0$
elsewhere.
By Lemma~\ref{L76}, $Y_{k}=o_{p}(N_k^{1/2})$.
Also, using \eqref{JW-bd}, one can conclude that
$\mathbf{u}_{0}^{\prime}W_{k}\mathbf{u}_{0}\geq{\sigma}_0^*+o_{p}(1)$
and hence,
$(1+t_{0}Y_{k})^{-1}t_{0}=O_{p}(N_k^{-1/2})$, a.s. ($\px$),
where ${\sigma}^*_0>0$ is the smallest eigenvalue of $ [{\sigma
}(\mathbf{0}
)]^2c_*^{-2}\Ga^*$.
Hence, it follows that $t_{0}=\llVert\b_{\theta_{0}}\rrVert
=O_{p}(N_k^{-1/2})$,
proving the claim.
\end{pf}

By the \hyperref[claim]{Claim} and Lemma~\ref{L76},
%
%
\begin{equation}
\quad\max_{1\leq j\leq N_k}|\g_{j}|\leq\llVert \b_{\theta
_{0}}
\rrVert Y_{k} =O_{p} \bigl(N_k^{-1/2}
\bigr)o_{p} \bigl(N_k^{1/2} \bigr)=o_{p}(1)
\qquad \mbox{a.s. ($\px$).} \label{gmax}
\end{equation}
Next, we obtain a stochastic approximation to $\b_{\theta_0}$. Using
$F_{k}(\theta_{0},\b_{\theta_{0}}) =0$, note that
\begin{eqnarray*}
0&=& N_k^{-1}\sum_{j=1}^{N_k}
\frac{G_{\theta_{0}}
(\w_{j})\tI_{k}(\w_{j})}{1+\b_{\theta_{0}}^{\prime} G_{\theta
_{0}}(\w_{j})
\tI_{n}(\w_{j})}
\\
&=& N_k^{-1}\sum
_{j=1}^{N_k}G_{\theta_{0}}(
\w_{j})\tI_{n}(\w _{j}) \biggl[1-
\g_{j} +\frac{\g_{j}^{2}}{1+\g_{j}} \biggr]
\\
&=& J_{k}-W_{k}\b_{\theta_{0}}+N_k^{-1}
\sum_{j=1}^{N_k} \frac{G_{\theta_{0}}(\w_{j})\tI_{k}(\w_{j})\g_{j}^{2}}{1+\g_{j}}.
\end{eqnarray*}
Therefore, we have the representation
%
%
\begin{equation}
\b_{\theta_{0}}=(W_{k})^{-1}J_{k}+
\eta_{k}, \label{beta}
\end{equation}
where, using condition (C.3), Lemma~\ref{L75}, the \hyperref[claim]{Claim},
and \eqref{gmax},
$
\llVert\eta_{k}\rrVert \leq$ $
Y_{k}\llVert\b_{\theta_{0}}\rrVert^{2} \llVert W_{k}\rrVert^{-1}
\{N_k^{-1}\sum_{j=1}^{N_k}\llVert G_{\theta_{0}}(\w_j)\rrVert^{2}
\tI_{k}^{2}(\w_{j})\}\{\max_{1\leq j\leq N_k}
(1 - |\g_{j}|)^{-1}\}
=o_{p}(N_k^{1/2}) O_{p}(N_k^{-1})O_{p}(1)O_{p}(1)O_{p}(1)=o_{p}(N_k^{-1/2})$,
a.s. ($\px$).
For $\llVert\b_{\theta_{0}}\rrVert Y_{k}<1$, applying a Taylor series
expansion, we have
\[
\log{(1+\g_{j})}=\g_{j}-\g_{j}^{2}/2+
\Delta_{j},
\]
where
$|\Delta_{j}|\leq\llVert\b_{\theta_{0}}\rrVert^{3}
Y_{k}\llVert G_{\theta_{0}}(\w_{j})\rrVert^{2}\tI_{k}^{2}(\w_{j})(1
-\llVert\b_{\theta_{0}}\rrVert Y_{k})^{-3}$
for all $1\leq j \leq N_k$.
Also, by Lemmas \ref{L73}--\ref{L75}, (C.3)
and \eqref{JW-bd},
$
N_kJ_{k}^{\prime}(W_{k})^{-1}J_{k} \stackrel{d}{\longrightarrow}
2\chi^{2}_{p}$ and
\begin{eqnarray*}
\sum_{j=1}^{N_k}|\Delta_{j}| &
\leq& N_k\llVert\b_{\theta_{0}}\rrVert^{3}Y_{k}
\bigl(1-\llVert\b _{\theta
_{0}}\rrVert Y_{k} \bigr)^{-3}
\Biggl\{N_k^{-1}\sum_{j=1}^{N_k}
\bigl\llVert G_{\theta
_{0}}(\w_{j}) \bigr\rrVert^{2}
\tI_{k}^{2}(\w_{j}) \Biggr\}
\\
&=&N_k O_{p}\bigl(N_k^{-3/2}
\bigr)o_{p}\bigl(N_k^{1/2}\bigr)O_{p}(1)O_{p}(1)=o_{p}(1),
\end{eqnarray*}
a.s. ($\px$).
Hence, it follows that\vspace*{2pt}
\begin{eqnarray*}
-\log\R_{n_k}(\theta_0) &\equiv& -\log\R_{k}(
\theta_0) = \sum_{j=1}^{N_k}
\log{(1+\g_{j})}
\\
&=& \Biggl[\sum_{j=1}^{N_k}
\g_{j}-2^{-1}\sum_{j=1}^{N_k}
\g _{j}^{2} \Biggr] +\sum_{j=1}^{N_k}
\Delta_{j}
\\
&=& \bigl[ \beta_{\theta_0}'[N_kJ_{k}]
- 2^{-1}N_k\beta_{\theta
_0}'W_{k}
\beta_{\theta_0} \bigr] +\sum_{j=1}^{N_k}
\Delta_{j}
\\
&=& 2^{-1}N_kJ_{k}^{\prime}(W_{k})^{-1}J_{k}
+ o_p(1) \stackrel{d} {\longrightarrow} \chi^2_p.
\end{eqnarray*}

This completes the proof of Theorem~\ref{thmm-1}.\vspace*{2pt}
\end{pf*}

\begin{pf*}{Proof of Theorem~\ref{thmm-2}}
By conditions on
$c_n$, $N$ and $\la_n$ in the MID case of part (a),\vspace*{2pt}
%
%
\begin{eqnarray}\label{mid-1}
b_n^2&\sim& Nc_n^{-2} \quad\mbox{and}
\nonumber
\\[-8pt]
\\[-8pt]
\nonumber
 \Biggl\|\Si_n - 2c_n^{-2} \bigl[{
\sigma}(\mathbf{0}) \bigr]^2\sum_{k=1}^NG_{\te
_0}(
\w_{kn}) G_{\te_0}(\w_{kn})' \Biggr\|& =& o
\bigl(b_n^2 \bigr),\vspace*{2pt}
\end{eqnarray}
where $c_n^{-1} = o(1)$. Thus, $b_n$ has a slower growth rate
in this case compared to the PID case. As in the proof of
Theorem~\ref{thmm-1}, it is enough to show that $-\log\R_{n_k}(\te_0)
\stackrel{d}{\longrightarrow} \chi^2_p$ through some subsequence $\{
n_k\}$ of a given
subsequence $\{n_i\}$. Indeed, the subsequence $\{n_k\}$
is extracted using (C.3)(iii) as before so that~\eqref{pid-2}
holds. Let $Y_k$, $\b_{\te_0}$
and $\g_j$ be as defined in the proof of Theorem~\ref{thmm-1}, and
here we continue to use the convention that
the subscript $n_k$ is replaced by~$k$, as before.
Next, redefine $J_k$ and $W_k$ as
$J_{k}= b_k^{-2}\sum_{j=1}^{N_k}G_{\theta_{0}}(\w_{kn})\tI_{k}(\w_{j})$
and
$
W_{k}
=
b_k^{-2}\sum_{j=1}^{N_k}G_{\theta_{0}}(\w_{j})G_{\theta_{0}}(\w
_{j})^{\prime}
\tI_{k}^{2}(\w_{j})$ where, following the convention,
we write $b_k =b_{n_k}$.
Then, by \eqref{pid-2}, \eqref{mid-1}, Lemma~\ref{L75}
and (C.4),\vspace*{2pt}
\begin{eqnarray*}
\bigl\|W_k - \bigl[{\sigma}(\mathbf{0}) \bigr]^2 \Ga^* \bigr\|
&=&o(1) \quad\mbox{and}
\nonumber
\\[-8pt]
\\[-8pt]
\nonumber
 b_{k}J_{k} &\cid& N \bigl(0, 2
\bigl[{ \sigma}(\mathbf{0}) \bigr]^2\Ga^* \bigr).
\end{eqnarray*}
Further, retracing the proof of Theorem~\ref{thmm-1}
and using Lemmas \ref{L75}--\ref{L77},
one can conclude that a.s. ($\px$),
$\|\b_{\te_0}\| = O_p(b_k^{-1})$ (cf. the \hyperref[claim]{Claim}),
$\max\{ |\g_j| \dvtx1\leq j \leq N_k\} = o_p(1)$ [cf. \eqref{gmax}]
and the representation \eqref{beta}\vadjust{\goodbreak} holds with
$\eta_k = o_p(b_k^{-1})$.
Hence, it follows that
\begin{eqnarray*}
-\log\R_{k}(\theta_0) 
&=& \Biggl[
\sum_{j=1}^{N_k}\g_{j}-2^{-1}
\sum_{j=1}^{N_k}\g _{j}^{2}
\Biggr] +\sum_{j=1}^{N_k}
\Delta_{j}
\\
&=& \bigl[ \beta_{\theta_0}'\bigl[b_k^2J_{k}
\bigr] - 2^{-1}b_k^2\beta_{\theta
_0}'W_{k}
\beta _{\theta_0} \bigr] +\sum_{j=1}^{N_k}
\Delta_{j}
\\
&=& 2^{-1}b_k^2J_{k}^{\prime}(W_{k})^{-1}J_{k}
+ o_p(1) 
\stackrel{d} {\longrightarrow}
\chi^2_p.
\end{eqnarray*}
%
%
This completes the proof of Theorem~\ref{thmm-2}(a).

Next, consider part (b). Note that in this MID case,
$N c_n^{-2} \ll\la_n^{\ka d}$
and hence, $b_n^2\sim\la_n^{\ka d}$. Also, using the boundedness of
$\|
G_{\te_0}(\cdot)\|$
over $\bbr^d$ and conditions (C.3)(i), (ii), (iv) and the DCT, one gets
\[
\bigl\|\Si_n - 2 \la_n^{\ka d} \Ga \bigr\| = o \bigl(
\la_n^{\ka d} \bigr),
\]
where $\Ga\equiv\int G_{\te_0}(\w) G_{\te_0}(\w)' K^2\phi^2 (\w
) \,d\w$
is nonsingular. Now retracing the proofs of Theorems \ref{thmm-1}
and \ref{thmm-2}(a) (with
$\{n_k\}$ replaced by the full sequence $\{n\}$), one can show that
$
-2\log\R_{n}(\theta_0)
= b_n^2 J_{0n}^{\prime}(W_{0n})^{-1}J_{0n}
+ o_p(1)$,
where
$J_{0n}= b_n^{-2}\sum_{j=1}^{N}G_{\theta_{0}}(\w_{jn})\tI_{n}(\w_{j})$
and
$
W_{0n}
=
b_n^{-2}\sum_{j=1}^{N}G_{\theta_{0}}(\w_{jn})G_{\theta_{0}}(\w
_{jn})^{\prime}
\tI_{n}^{2}(\w_{jn})$. Note that by Lemma~\ref{L75} and the
fact that $b_n^2\sim\la_n^{\ka d}$, we have
\[
\bigl\|W_n - b_n^{-2}\Si_n \bigr\| = o(1),
\]
which is different from the
previous two cases covered by Theorems \ref{thmm-1}
and~\ref{thmm-2}(a) [where $\|W_n - 2^{-1} (b_n^{-2}\Si_n)\|
= o(1)$]. In view of (C.4), this implies that
$b_n^2 J_{0n}^{\prime} (W_{0n})^{-1} J_{0n}
\stackrel{d}{\longrightarrow} \chi^2_p$,
proving part (b).
\end{pf*}

\begin{remark}
\label{re7.1} From the proof of Theorems \ref
{thmm-1}--\ref{thmm-2},
it follows that
the different scalings in the two cases are required by the dominant
term in the asymptotic variance of the sum $\sum_{k=1}^N G_{\te_0}(\w
_{kn}) I_n^*(\w_{kn})$ and the automatic variance stabilization factor,
both of which arise from the inner mechanics of the SFDEL. Under PID
and under ``slow'' MID, the leading term is given by $Nc_n^{-1}{\sigma
}(\mathbf{0}
)$, which is of a larger order of magnitude than $\la_n^{\ka d}$. When
the infilling rate is high, that is, $Nc_n^{-2}\ll\la_n^{\ka d}$, the
other term involving the spectral density of the $Z(\cdot)$-process
dominates (as in the case of regularly spaced time series FDEL) and the
standard scaling by $-2$ is appropriate.\
\end{remark}
%

\begin{pf*}{Proof of Theorem~\ref{thmm-4}} By Lemma~\ref{L75} and the fact
that $\hat{{\sigma}}_n(\mathbf{0}) - {\sigma}(\mathbf{0})
= O_p(\la_n^{-d/2})$ (cf. Lahiri \cite{lahiri2003a}), under the
conditions of Theorem~\ref{thmm-2}(b),
\begin{eqnarray*}
a_n(\te_0) &=& \frac{b_n^{-2} \sum_{j=1}^N \|G_{\te_0} (\w_{jn})\|
^2 \tI
_n^2(\w_{jn})}{
b_n^{-2} \sum_{j=1}^N \|G_{\te_0} (\w_{jn})\|^2 I_n^2(\w_{jn})}
\\
&=& \frac{b_n^{-2} \sum_{j=1}^N \|G_{\te_0} (\w_{jn})\|^2 [2K^2\phi
^2(\w
_{jn})] + o_p(1)}{
b_n^{-2} \sum_{j=1}^N \|G_{\te_0} (\w_{jn})\|^2 [2K^2\phi^2(\w_{jn})]
+ o_p(1)
}
= 1+o_p(1)
\end{eqnarray*}
while under the conditions of Theorems \ref{thmm-1}
and \ref{thmm-2},
\begin{eqnarray*}
a_{n_k}(\te_0) &=& \frac{b_k^{-2}
\sum_{j=1}^{N_k} \|G_{\te_0} (\w_{k})\|^2 [c_{n_k}^{-2} [{\sigma
}(\mathbf{0}
)]^2] + o_p(1)}{
b_k^{-2}
\sum_{j=1}^{N_k} \|G_{\te_0} (\w_{k})\|^2 [2c_{n_k}^{-2} [{\sigma
}(\mathbf{0}
)]^2] + o_p(1)}
\\
&=& 2^{-1} \bigl(1+o_p(1) \bigr).
\end{eqnarray*}

Now combining this with the proofs of Theorems \ref{thmm-1}--\ref
{thmm-2}, one
can complete the
proof of Theorem~\ref{thmm-4}.
\end{pf*}

\section*{Acknowledgements} The authors are grateful to three reviewers
and an Associate Editor for thoughtful comments
and constructive criticism which led to significant improvements in the
manuscript.

\begin{supplement}
\stitle{Supplement to ``A frequency domain empirical likelihood
method for irregularly spaced spatial data''}
\slink[doi]{10.1214/14-AOS1291SUPP} 
\sdatatype{.pdf}
\sfilename{aos1291\_supp.pdf}
\sdescription{Details of proofs and additional simulation results.}
\end{supplement}


%



\printaddresses
\end{document}